\documentclass{gtart}

\def\ifplaintex{\expandafter\ifx\csname documentclass\endcsname\relax}

\def\gtp{{\mathsurround=0pt\it $\cal G\mskip-2mu$eometry \&\ 
$\cal T\!\!$opology $\cal P\!$ublications}}  

\def\recd{{\small Received:\qua\receiveddate\ifx\reviseddate\relax
\else\qquad Revised:\qua\reviseddate\fi\par}} 

\def\lognumber#1{\def\thelognumber{#1}}
\def\volumenumber#1{\def\thevolumenumber{#1}}
\def\volumeyear#1{\def\thevolumeyear{#1}}
\def\papernumber#1{\def\thepapernumber{#1}}
\def\pagenumbers#1#2{\def\startpage{#1}\def\finishpage{#2}}
\def\published#1{\def\publishdate{#1}}

\def\received#1{\def\receiveddate{#1}}
\def\revised#1{\def\reviseddate{#1}}
\def\accepted#1{\def\accepteddate{#1}}

\let\\\par\let\thelognumber\relax\let\thevolumenumber\relax
\let\thepapernumber\relax\let\thevolumeyear\relax\let\startpage\relax
\let\finishpage\relax\let\publishdate\relax\let\receiveddate\relax
\let\reviseddate\relax\let\accepteddate\relax\let\theasciititle\relax
\let\theasciiauthors\relax
\let\theasciiabstract\relax

\let\theasciiemail\relax

\ifplaintex
\font\logobig=cmssbx10 scaled 3836
\font\logomed=cmssbx10 scaled 2557
\else
\font\logobig=cmssbx10 scaled 4200
\font\logomed=cmssbx10 scaled 2800
\fi

\long\def\makeagttitle{   
\count0=\startpage
\agt\hfill      
\hbox to 45truept{\vbox to 0pt{\vglue -13truept{\logomed A\kern -.37em{\logobig 
T}\kern -.38em G}\vss}\hss}
\break
{\small Volume \thevolumenumber\ (\thevolumeyear)
\startpage--\finishpage\nl
Published: \publishdate}

\vglue .25truein

{\parskip=0pt\leftskip 0pt plus
1fil\def\\{\par\smallskip}{\Large\bf\thetitle}\par\medskip} \vglue
0.05truein

{\parskip=0pt\leftskip 0pt plus 1fil\def\\{\par}{\sc\theauthors}
\par\medskip}%
 
\vglue 0.03truein 

{\small\leftskip 25truept\rightskip 25truept{\bf Abstract}\stdspace\theabstract

{\bf AMS Classification}\stdspace\theprimaryclass
\ifx\thesecondaryclass\relax\else; \thesecondaryclass\fi\par
{\bf Keywords}\stdspace \thekeywords\par}\vglue 7truept

}   

\ifplaintex
\font\phead=cmsl9 scaled 950
\font\pnum=cmbx10 scaled 913
\font\pfoot=cmsl9 scaled 950
\headline{\vbox to 0pt{\vskip -4.5mm\line{\small\phead\ifnum
\count0=\startpage ISSN 1472-2739 (on-line) 1472-2747 (printed)
\hfill {\pnum\folio}\else\ifodd\count0\def\\{ }%
\ifx\theshorttitle\relax\thetitle\else\theshorttitle\fi\hfill{\pnum\folio}
\else\def\\{ and }{\pnum\folio}\hfill\ifx\theshortauthors\relax\theauthors
\else\theshortauthors\fi\fi\fi}\vss}}
\footline{\vbox to 0pt{\vglue 0mm\line{\small\pfoot\ifnum\count0=\startpage
\copyright\ \gtp\hfill\else
\agt, Volume \thevolumenumber\ (\thevolumeyear)\hfill\fi}\vss}}
\else
\font=cmsl9 scaled 1050
\font\lnum=cmbx10 
\font=cmsl9 scaled 1050
\makeatletter
\def\@oddhead{{\small\lhead\ifnum\count0=\startpage ISSN 1472-2739 
(on-line) 1472-2747 (printed)\hfill {\lnum\number\count0}\else\ifodd\count0
\def\\{ }\ifx\theshorttitle\relax \thetitle \else\theshorttitle\fi\hfill
{\lnum\number\count0}\else\def\\{ and }{\lnum\number\count0}
\hfill\ifx\theshortauthors\relax 
\theauthors\else\theshortauthors\fi\fi\fi}}\def\@evenhead{\@oddhead}
\def\@oddfoot{\small\lfoot\ifnum\count0=\startpage\copyright\ \gtp\hfill\else
\agt, Volume \thevolumenumber\ (\thevolumeyear)\hfill\fi}
\def\@evenfoot{\@oddfoot}
\makeatother
\fi
\let\maketitlepage\makeagttitle
\let\makeshorttitle\maketitlepage
\let\maketitle\maketitlepage

\newwrite\gtoutfile
\long\gdef\makeheadfile{  
{\def\\{, }\def\s{ }
\immediate\openout\gtoutfile head.xxx
\immediate\write\gtoutfile{To: math@arxiv.org}
\immediate\write\gtoutfile{Subject: put OR rep NNNNN:ppppp}
\immediate\write\gtoutfile{--text follows this line--}
\immediate\write\gtoutfile{Proxy-for: \ifx\theasciiauthors\relax
\theauthors\else\theasciiauthors\fi\s<\ifx\theasciiemail\relax\theemail\else\theasciiemail\fi>}
\immediate\write\gtoutfile{\noexpand\\}
\immediate\write\gtoutfile{Authors: \ifx\theasciiauthors\relax
\theauthors\else\theasciiauthors\fi}
{\def\\{ }\immediate\write\gtoutfile{Title: \ifx\theasciititle\relax
\thetitle\else\theasciititle\fi}}
\immediate\write\gtoutfile{Subj-class: GT or SG, GR etc}
\immediate\write\gtoutfile{MSC-class: \theprimaryclass\ifx\thesecondaryclass\relax\else, \thesecondaryclass\fi}
\immediate\write\gtoutfile{Journal-ref: Algebr. Geom. Topol. \thevolumenumber\s
(\thevolumeyear) \startpage-\finishpage}
\immediate\write\gtoutfile{Comments: Published by Algebraic and
Geometric Topology at}
\immediate\write\gtoutfile{\s\s\s  http://www.maths.warwick.ac.uk/agt/AGTVol\thevolumenumber/agt-\thevolumenumber-\thepapernumber.abs.html}
\immediate\write\gtoutfile{\noexpand\\}
\immediate\write\gtoutfile{}
\ifx\theasciiabstract\relax
\immediate\write\gtoutfile{\theabstract}\else
\immediate\write\gtoutfile{\theasciiabstract}\fi
\immediate\write\gtoutfile{}
\immediate\write\gtoutfile{\noexpand\\}
\immediate\write\gtoutfile{}
\immediate\closeout\gtoutfile}}  

\def\maketitlepage{\makeagttitle\makeheadfile}
\let\makeshorttitle\maketitlepage
\let\maketitle\maketitlepage

\lognumber{1}
\volumenumber{3}
\volumeyear{2003}
\papernumber{1}
\published{24 January 2002}
\pagenumbers{1}{31}
\received{8 October 2002}
\revised{22 December 2002}
\accepted{10 January 2003}

\usepackage{amssymb,amsmath,graphics}

\newtheorem{thm}{Theorem}[section]  
\newtheorem{lem}[thm]{Lemma}        
\newtheorem{cor}[thm]{Corollary}
\theoremstyle{definition}
\newtheorem{defn}[thm]{Definition}  
\newtheorem*{rem}{Remark}           

\newcommand{\thmref}[1]{Theorem~\ref{#1}}
\newcommand{\secref}[1]{\S\ref{#1}}

\newcommand{\wt}{\widetilde}
\newfont{\cyr}{wncyr10}
\newcommand{\Lob}{\mbox{\cyr L}}

\renewcommand{\L}{{\mathcal L}}

\newcommand{\cfig}[3]{
\begin{figure}[ht!]
\center{\scalebox{.7}{\includegraphics{#1}}}
\caption{#3}
\label{#2}
\end{figure}
}
\newcommand{\afig}[3]{
\begin{figure}[ht!]
\center{\scalebox{.6}{\includegraphics{#1}}}
\caption{#3}
\label{#2}
\end{figure}
}
\newcommand{\bfig}[3]{
\begin{figure}[ht!]
\center{\scalebox{.8}{\includegraphics{#1}}}
\caption{#3}
\label{#2}
\end{figure}
}

\newcommand{\figref}[1]{Figure~\ref{#1}}

\begin{document}

\title{The Regge symmetry is a scissors congruence\\in hyperbolic space}
\authors{Yana Mohanty}
\address{Department of Mathematics, University of California at San 
Diego\\La Jolla, CA 92093-0112, USA }

\email{mohanty@math.ucsd.edu}
\begin{abstract} 
We give a constructive proof that the Regge symmetry is a scissors congruence in hyperbolic space.
The main tool is Leibon's construction for computing the volume of a general
hyperbolic tetrahedron. The proof consists of identifying the key elements in
Leibon's construction and permuting them.
\end{abstract}


\maketitle

\section{Introduction}\label{s:intro}
The Regge symmetries are a family of involutive linear transformations on the six edges of
a tetrahedron. They are defined as follows.
\begin{defn}\label{d:regge}
Let $T(A,B,C,A^\prime,B^\prime,C^\prime)$ denote a tetrahedron as shown in Figure~\ref{f:generic_tet}.
Define
\begin{equation}\label{e:R_a}
R_a(T(A,B,C,A^\prime,B^\prime,C^\prime))=T(A,s_a-B,s_a-C,A^\prime,s_a-B^\prime,s_a-C^\prime),
\end{equation}
where
$s_a=(B+C+B^\prime+C^\prime)/2$. Similarly, define
\begin{equation}\label{e:R_b}
R_b(T(A,B,C,A^\prime,B^\prime,C^\prime))=T(s_b-A,B,s_b-C,s_b-A^\prime,B^\prime,s_b-C^
\prime),
\end{equation}
and
\begin{equation}\label{e:R_c}
R_c(T(A,B,C,A^\prime,B^\prime,C^\prime))=T(s_c-A,s_c-B,C,s_c-A^\prime,s_c-B^\prime,C^\prime),
\end{equation}
where $s_b=(A+C+A^\prime+C^\prime)/2$ and
$s_c=(A+B+A^\prime+B^\prime)/2$. Then $R_a$, $R_b$, and $R_c$ generate the family of {\em
Regge symmetries} of $T(A,B,C,A^\prime,B^\prime,C^\prime)$. \end{defn}
\bfig{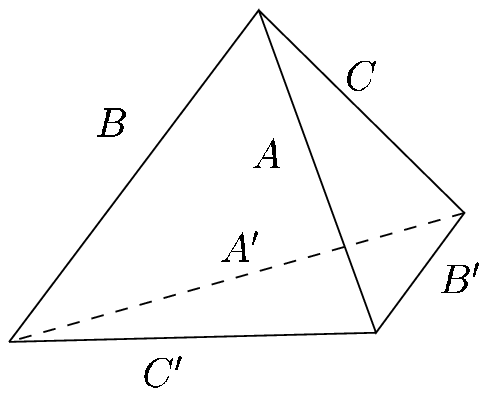}{f:generic_tet}{Tetrahedron
$T(A,B,C,A^\prime,B^\prime,C^\prime)$ with its dihedral angles denoted by
letters} Any two maps out of $R_a$, $R_b$, and $R_c$, together with the 
tetrahedral symmetries, form a group isomorphic to $S_3\times S_4$
\cite{Roberts}.

The Regge symmetries first arose in conjunction
with the $6j$-symbol, which is a real number that can be associated to a labeling of the
six edges of a tetrahedron by irreducible representations of $SU(2)$. In the 1960's,
Tullio Regge discovered that the $6j$-symbols are invariant under the linear
transformations generated by $R_a$, $R_b$, and $R_c$. Expanding on his work in
\cite{Roberts}, Justin Roberts explored the effect of the Regge symmetries on Euclidean
tetrahedra associated with the $6j$-symbols. He found that the volumes of Euclidean
tetrahedra as well as their Dehn invariants remain unchanged under the action of the Regge
symmetries. Therefore, the Regge symmetries give rise to a family of scissors congruent
Euclidean tetrahedra. In the hyperbolic case, it was unknown until now whether the Regge
symmetries preserve equidecomposability, since the conjecture concerning the completeness
of the volume and Dehn invariant as scissors congruence invariants is still open. In this
paper, we show that the Regge symmetries do indeed generate a family of scissors
congruent tetrahedra by an explicit construction. The construction is based on a volume
formula for a hyperbolic tetrahedron first developed by Jun Murakami and
Masakazu Yano in \cite{Murakami}, and later geometrically interpreted and generalized by
Gregory Leibon in \cite{Leibon}.

\subsection{Outline of the argument}
Given a finite hyperbolic tetrahedron $T$, Leibon's formulas give a geometric
decomposition of $2T$ (2 copies of $T$) of the form
$\coprod_{j=1}^{16}\L(\theta_j)$, where $\L(\theta)$ is defined to be
half of an isosceles (bilaterally symmetric) ideal tetrahedron with apex angle
$2\theta$, as shown in \figref{f:isosceles}.
{\begin{figure}
\center{\includegraphics{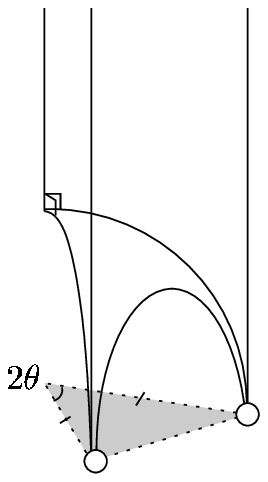}}
\vspace{-30pt}
\caption{Bilaterally symmetric 3/4-ideal tetrahedron $\L(\theta)$ in the half-space model.
The shaded isosceles triangle is the shadow cast by the tetrahedron onto the
plane at infinity. The circles represent vertices at infinity.}
\label{f:isosceles} \end{figure}
}
Leibon's construction is based on
the idea of extending all the edges of the tetrahedron to infinity and dissecting
the resulting polyhedron into 6 ideal tetrahedra and an ideal octahedron. The
construction is identical for hyperideal tetrahedra, as shown in
\secref{s:warm_up}, and in this case the convex hull of the 12 ideal vertices is
combinatorially equivalent to a truncated tetrahedron. Such a polyhedron can be
triangulated by tetrahedra that do not intersect each other, and for this reason
it is much easier to see the scissors congruence proof in the hyperideal case.
It turns out that there is a ``dual'' dissection which
results in an octahedron whose dihedral angles are supplementary to the angles
of the original octahedron. Moreover, it will be shown in
\secref{s:gregs_formula} that $2T$ can be constructed just from the original
octahedron and its dual. In order to get to a decomposition of $T$, we use
Dupont's result in \cite{Dupont_book} that the group of hyperbolic polyhedra is uniquely
2-divisible. This allows us to literally halve Leibon's construction by slicing
through each of the $\L(\theta_j)$ along its plane of symmetry. This gives the
decomposition $T=\coprod_{j=1}^{16} \L_h(\theta_j)$,
where $\L_h(\theta_j)$ is one of the bilaterally symmetric halves of
$\L(\theta)$. In \secref{s:regge_proof} we show how four of the
$\L_h(\theta_j)$ and be permuted so that the result is $R_b(T)$, the image
of $T$ under one of the Regge symmetries.

In an effort to make this paper self-contained, we have summarized the results
we will need from \cite{Leibon} in
\secref{s:warm_up} and \secref{s:gregs_formula}.
The exposition of these results in \cite{Leibon} is more
general, while the presentation here is geared specifically for what we will need to
demonstrate the
scissors congruence proof in \secref{s:regge_proof}.

I would like to thank Greg Leibon and Peter Doyle for generously sharing their
ideas with me. Most of all, thanks to Justin Roberts for introducing us all to
this subject.

\section{Warm-up for the development Leibon's set of formulas
for the volume of a hyperbolic tetrahedron}
\label{s:warm_up}
In this section we develop the basic idea that is central to Leibon's geometrization of
Murakami and Yano's formula. Let
$T_3(A,B,C)$ denote a 3/4-ideal tetrahedron with dihedral angles $A$, $B$, and $C$
 at its finite vertex.
We now extend the three edges meeting at the non-ideal vertex to infinity obtaining the
polyhedron $D$ shown in Figure~\ref{f:3_4ideal}. For simplicity, all the
hyperbolic polyhedra from now on will be shown in the Klein model, unless
specified otherwise.

\afig{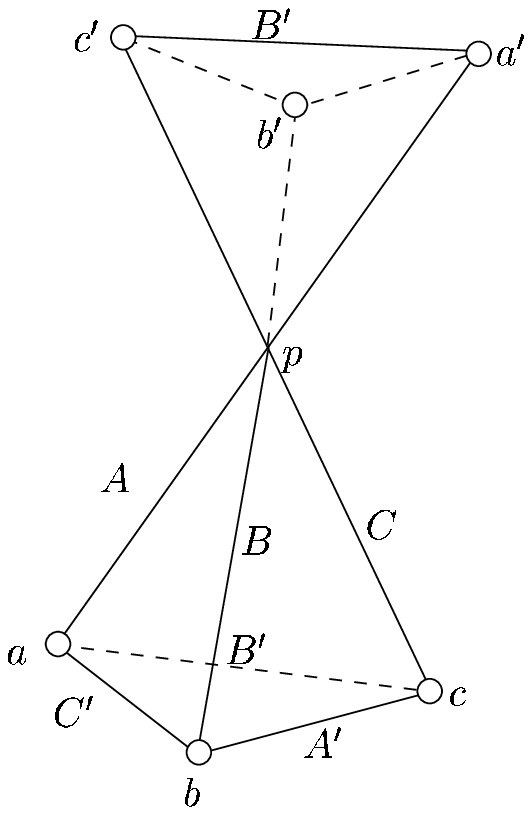}{f:3_4ideal}{Polyhedron $D$}

Notice that $D$ is symmetric about the point $p$, so that its volume
is twice that of $T_3(A,B,C)$. It will be convenient to view $D$ as a
simplicial complex which can be triangulated by oriented simplices as
follows.
\begin{equation}\label{e:3_ideal_a}
D=
\{a,b,c,c^{\prime}\}+\{a,a^{\prime},b^{\prime},c^{\prime}\}+\{a,b^{\prime} ,b,
c^{\prime}\}, \end{equation}
where $\{x_1,x_2,x_2,x_4\}$ denotes the oriented hyperbolic tetrahedron with
vertices $x_1$, $x_2$, $x_3$, and $x_4$.
\afig{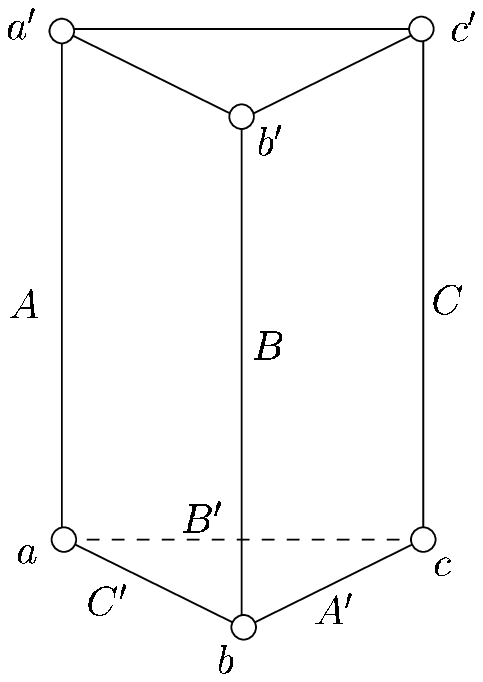}{f:prism}{An ideal hyperbolic prism}
Now consider the ideal hyperbolic prism, $P$, shown in Figure~\ref{f:prism}. Viewed as a
simplicial complex, $P$ can be triangulated as
\begin{equation}\label{e:prism}
P=
\{a,b,c,c^{\prime}\}+\{a,a^{\prime},b^{\prime},c^{\prime}\}+\{a,b^{\prime} ,b,
c^{\prime}\}. \end{equation}
Equations \eqref{e:3_ideal_a} and \eqref{e:prism} indicate that $P$ and $D$ are the same
object from the point of view of homology. Clearly, $P$ and $D$ are embedded in space
differently. In particular, $P$ is convex while $D$ is not. This is accounted for by the
term $\{a,b^{\prime},b,c^{\prime}\}$ in equations \eqref{e:3_ideal_a} and
\eqref{e:prism}. In the case of $D$,  $\{a,b^{\prime} ,b,c^{\prime}\}$ represents a
simplex with negative volume, while in the case of $P$, it represents a simplex with
positive volume.

Just as $D$ and $P$ are the same object when viewed as simplicial
complexes, they are also the same type of object from the point of view of hyperbolic
geometry. Both $D$ and $P$ are completely determined by the dihedral angles $A$, $B$, and
$C$. In
the case of $D$, $A+B+C>\pi$, while in the case of $P$, $A+B+C<\pi$. $P$ can be obtained
by a continuous deformation of $D$ which involves moving the point $p$ outside the sphere
at infinity (this is easiest seen in the Klein or the hyperboloid model). In this process
 the angles $A$, $B$, and $C$
decrease. Since the volume of $D$ and $P$ depends only on these three angles, by analytic
continuation $D$ and $P$ have the same volume formula. It is easier to see the
triangulation of $P$ since the three tetrahedra involved do not
intersect each other.

Using the fact that the opposite dihedral angles of an ideal hyperbolic
tetrahedron are equal, we find that the tetrahedra in
\figref{f:prism_split} are \begin{eqnarray}
\label{e:3_ideal_b}
\{a,b,c,c^{\prime}\}&=&T(A^\prime,B^\prime,C)\nonumber\\
\{a,a^{\prime},b^{\prime},c^{\prime}\}&=&T(A,B^\prime,C^\prime)\\
\{a,b,b^{\prime} ,c^{\prime}\}& =&T(C^\prime-C,B,\pi-B^\prime)\nonumber,
\end{eqnarray}
where $T(A,B,C)$ denotes an ideal hyperbolic tetrahedron with dihedral angles $A$,
$B$, and $C$.
\afig{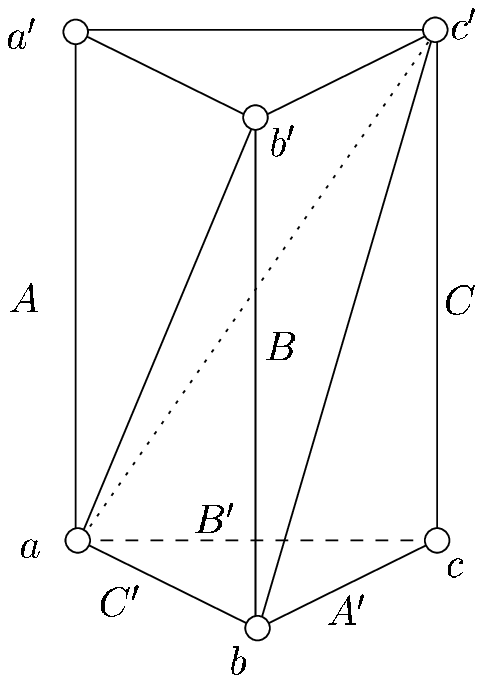}{f:prism_split}{A triangulation of an ideal hyperbolic
prism} Applying the condition that the sum of the dihedral angles at an ideal 
vertex is $\pi$, we obtain \begin{eqnarray}\label{e:3_ideal_others}
A^{\prime}=\frac{\pi+A-B-C}{2}\nonumber\\
B^{\prime}=\frac{\pi+B-A-C}{2}\\
C^{\prime}=\frac{\pi+C-A-B}{2}\nonumber.
\end{eqnarray}
By \eqref{e:prism}, \eqref{e:3_ideal_b} and the famous formula from
\cite{Milnor},
 \begin{equation}\label{e:allideal}
 V(T(\alpha,\beta,\gamma))=\Lob(\alpha)+\Lob(\beta)+\Lob(\gamma),
 \end{equation}
 where $\Lob(\theta):=-\int_{0}^{\theta}\log 2|\sin u|du$ is the {\em Lobachevsky function},
 we have
\begin{multline}\label{e:3_ideal_c}
2V(\{a,b,c,a^\prime,b^\prime,c^\prime\})=\Lob(A)+
\Lob(A^{\prime})+\Lob(B)+\Lob(B^{\prime})+
\Lob(C)+\Lob(C^{\prime})-\\ \Lob(\frac{\pi+A+B+C}{2}),
\end{multline}
where $\{a,b,c,a^\prime,b^\prime,c^\prime\}$ is either the non-convex prism of
Figure~\ref{f:3_4ideal} or the convex prism of Figure~\ref{f:prism}.

\section{Leibon's formulas for the volume of a hyperbolic tetrahedron}
\label{s:gregs_formula}
\subsection{The basic setup}\label{s:basic_setup}
We now extend the ideas developed in \secref{s:warm_up} to a hyperbolic tetrahedron
$T$ with finite vertices. We start by extending all the edges of $T$ to infinity. The
resulting polyhedron, $C$, is shown in Figure~\ref{f:wrapped_candy}, where
$T=\{p_1,p_2,p_3,p_4\}$. Clearly,
\begin{multline}\label{e:wrapped_candy}
V(T)=V(C)-\\
[V(\{c_1^\prime,a_1,b_1^\prime,p_1\})+V(\{a_1^\prime,b_2^\prime,c_1,p_2\})+
V(\{a_2^\prime,b_1,c_2^\prime,p_3\})+V(\{a_2,b_2,c_2,p_4\})] \end{multline}
Since the 4 tetrahedra on the right hand side of \eqref{e:wrapped_candy} are
all 3/4-ideal, their volume is given by \eqref{e:3_ideal_c}. Therefore, the main
task at hand is to calculate the volume of $C$. In order to triangulate $C$, we 
first note that it is the same object as $U$ (see
Figure~\ref{f:unwrapped_candy}) from the point of view of homology, following
the method of \secref{s:warm_up}. That is, $U$ can be obtained from $C$ by
pulling the points $p_1$, $p_2$, $p_3$, and $p_4$ outside the sphere at
infinity. Under this deformation the non-convex prisms
$\{c_1^\prime,a_1,b_1^\prime,c_2^\prime,a_2,b_2^ \prime\}$,
$\{a_1^\prime,b_2^\prime,c_1,a_2^\prime,b_1^\prime ,c_2\}$,
$\{a_2^\prime,b_1,c_2^\prime,a_1^\prime,b_2,c_1^ \prime\}$, and
$\{a_2,b_2,c_2,a_1,b_1,c_1\}$ become convex, but their volume formula stays the
same by analytic continuation. Thus
\cfig{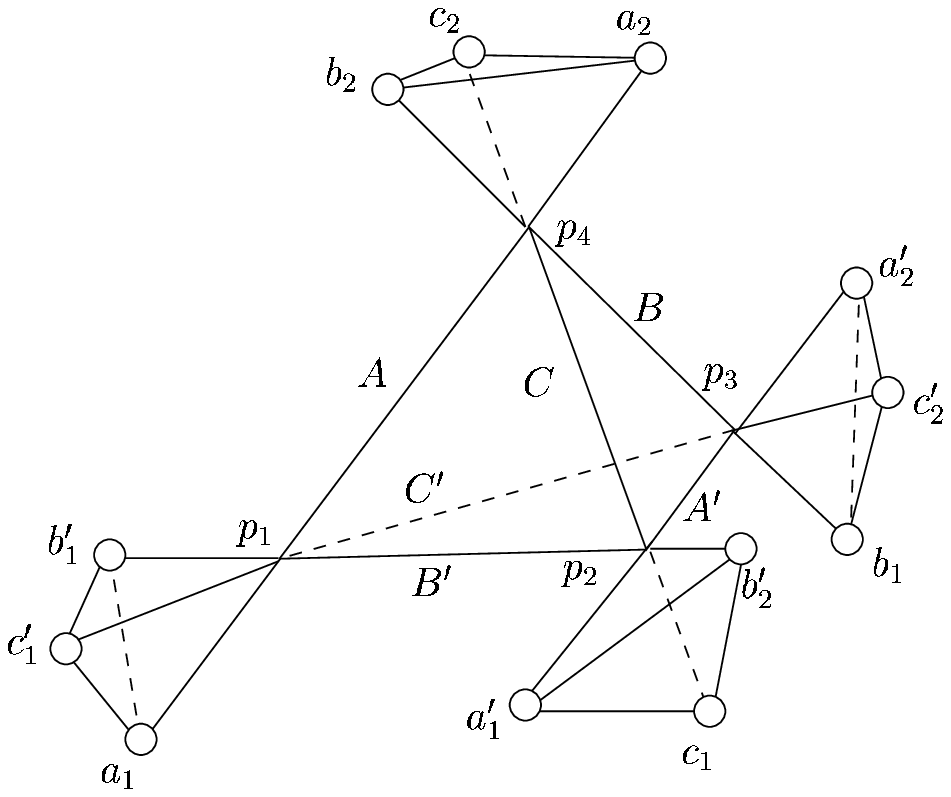}{f:wrapped_candy}{Polyhedron $C$}
\afig{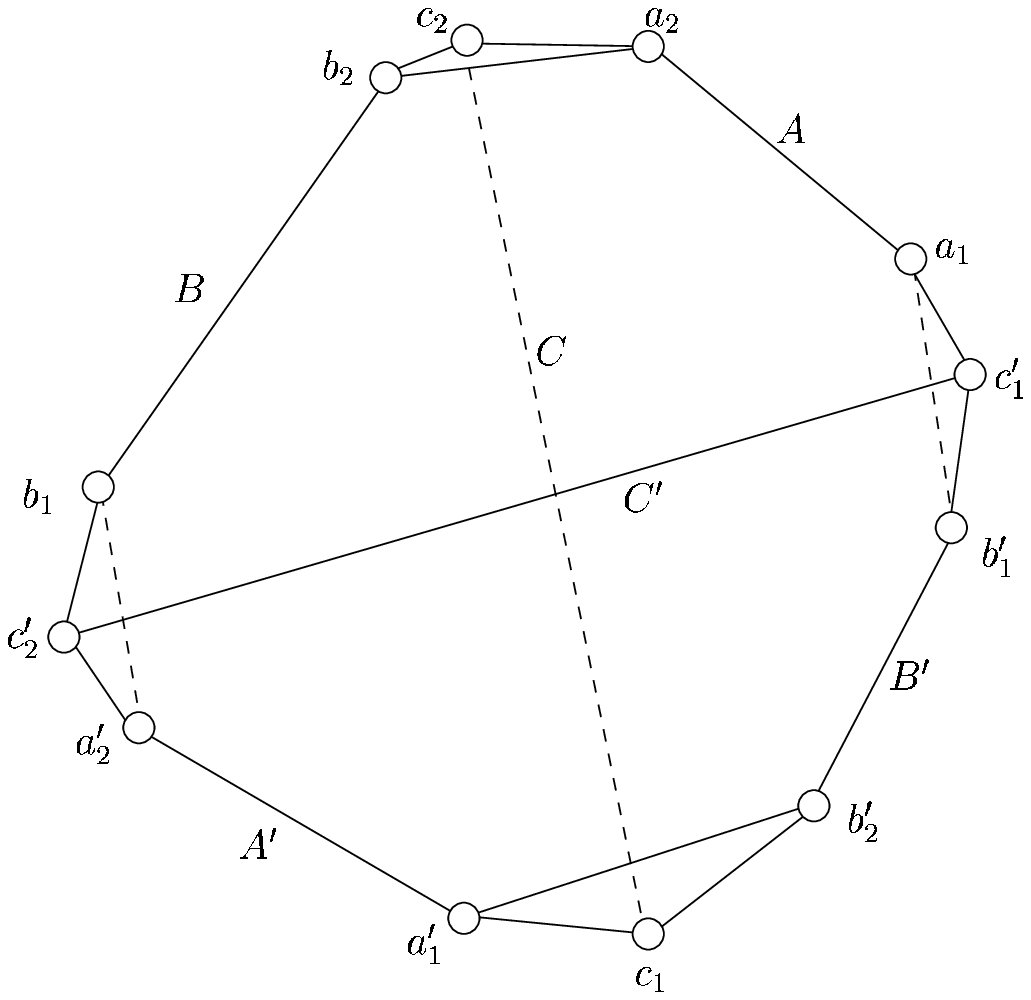}{f:unwrapped_candy}{Polyhedron $U$}
\begin{multline}\label{e:unwrapped_candy}
V(T)=V(U)-V(\{c_1^\prime,a_1,b_1^\prime,c_2^\prime,a_2,b_2^ \prime\})/2-
V(\{a_1^\prime,b_2^\prime,c_1,a_2^\prime,b_1^\prime ,c_2\})/2\\
-V(\{a_2^\prime,b_1,c_2^\prime,a_1^\prime,b_2,c_1^
\prime\})/2-V(\{a_2,b_2,c_2,a_1,b_1,c_1\})/2.
\end{multline}
Since a triangulation of $C$ by oriented simplices is also a triangulation of
$U$, we may as well work with $U$. This is easier than working with $C$ because
the tetrahedra in the triangulation of $U$ do not intersect each other.
Everything that follows applies equally well to finite as well as hyperideal
tetrahedra.

In order to compute the volume of $U$ we triangulate it using
only ideal tetrahedra  and the compute the volume of each of these tetrahedra
using \eqref{e:allideal}. It turns out that no matter which way $U$ is
triangulated, some of the tetrahedra comprising it will have dihedral angles
that are affine functions of the dihedral angles of $T$, while others will not.
Leibon has looked at a particular family of $2^6$ triangulations of $U$, each of
which divides up $U$ into six tetrahedra and one octahedron. The number $2^6$ comes
from the fact that the vertices of the octahedron can be chosen by selecting
exactly one of the vertices in each pair $\{a_1,a_2\}$, $\{b_1,b_2\}$,
$\{c_1,c_2\}$, $\{a_1^\prime,a_2^\prime\}$, $\{b_1^\prime,b_2^\prime\}$,
$\{c_1^\prime,c_2^\prime\}$. If one chooses both or neither of the vertices in
any pair, then it is impossible to form a non-degenerate octahedron with the
remaining vertices.

We will illustrate the computation of the volume of $U$ using the partial
triangulation shown in \figref{f:U_decomp_1}. The prism and three tetrahedra
resulting from this decomposition are shown in Figure~\ref{f:U_decomp_tet_lab},
while the octahedron $O$ is shown in Figure~\ref{f:octahedron}. The prism in Figure~\ref{f:U_decomp_tet_lab} can be decomposed into
three ideal tetrahedra, as demonstrated in \secref{s:warm_up}.
In Figure~\ref{f:U_decomp_tet_lab}, the dihedral angles other than those belonging to the
original tetrahedron $T$ were computed with the help of
\eqref{e:3_ideal_others}. For
example, the dihedral angle at the edge $\{a_1,c_1^\prime\}$ was computed by
considering the prism $\{c_1^\prime,a_1,b_1^\prime,c_2^\prime,a_2,b_2^
\prime\}$. Since the dihedral angles at each pair
of opposite edges of an ideal tetrahedron are equal, the labels in
Figure~\ref{f:U_decomp_tet_lab} are enough to completely determine each tetrahedron. The
dihedral angles of the octahedron $O$ then follow immediately.
\afig{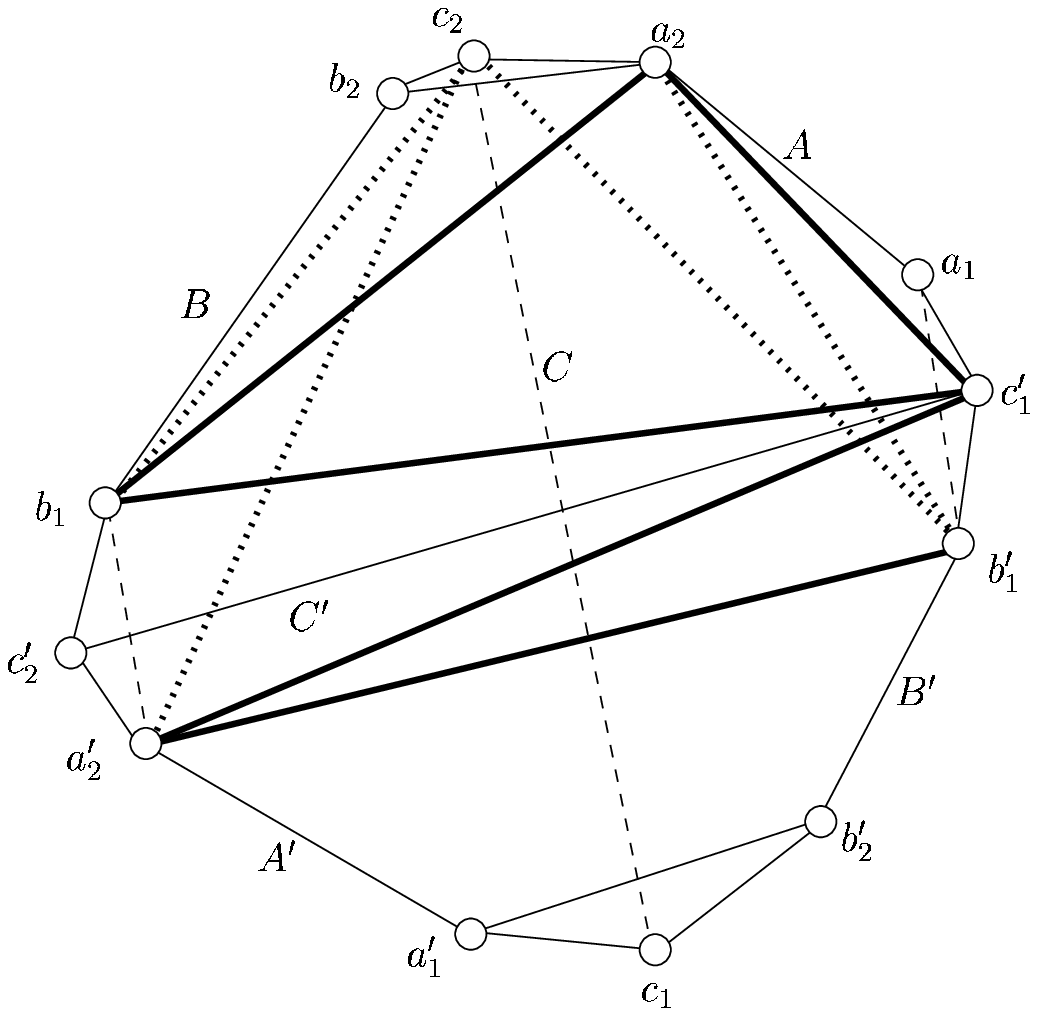}{f:U_decomp_1}{Decomposition of
polyhedron $U$. The edges of $U$ are shown with thin and dashed lines, while the
cuts of the decomposition are shown with thicker and dotted lines.}
\afig{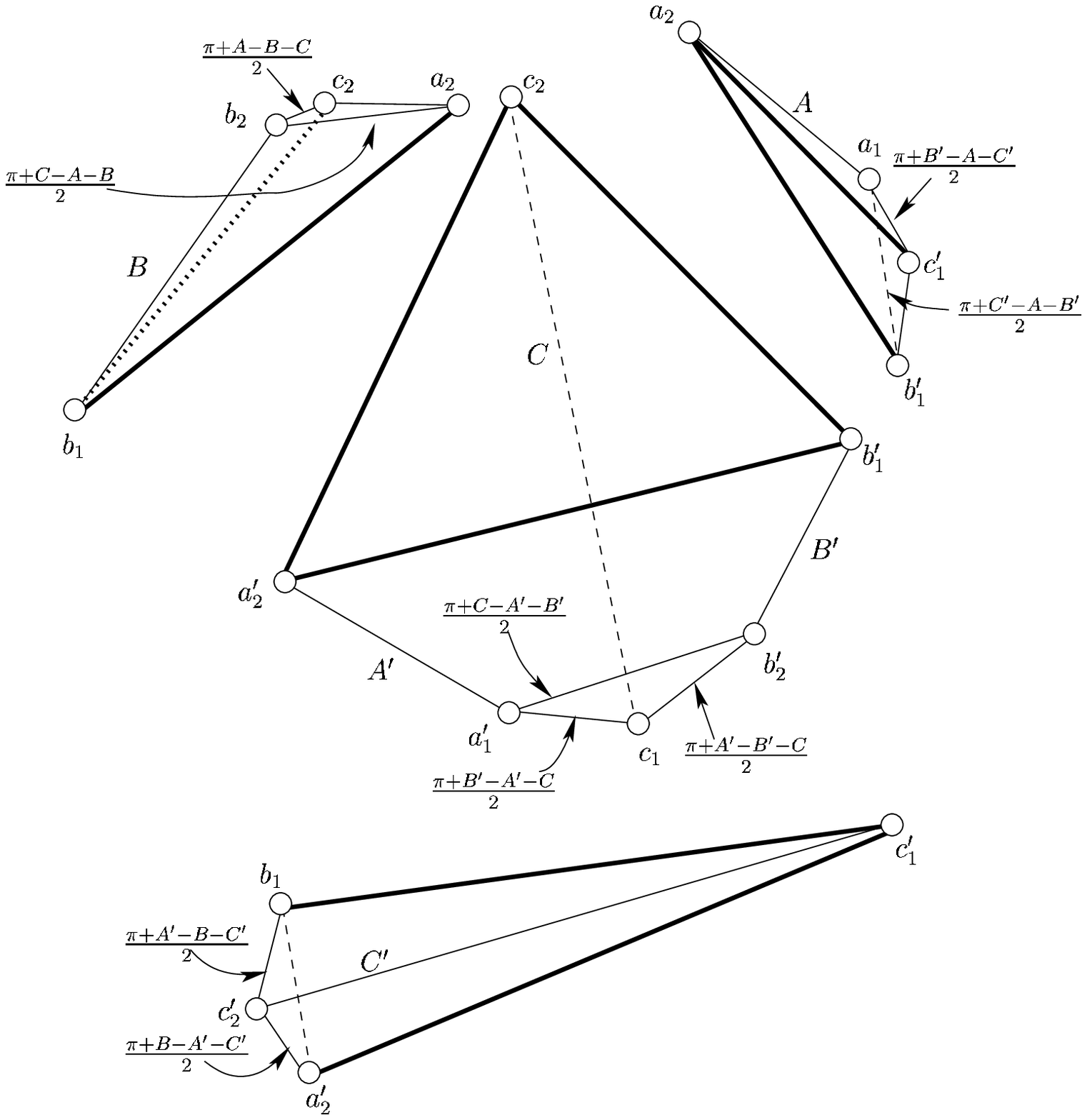}{f:U_decomp_tet_lab}{Decomposition of polyhedron
$U$, Part I: the tetrahedra}
\cfig{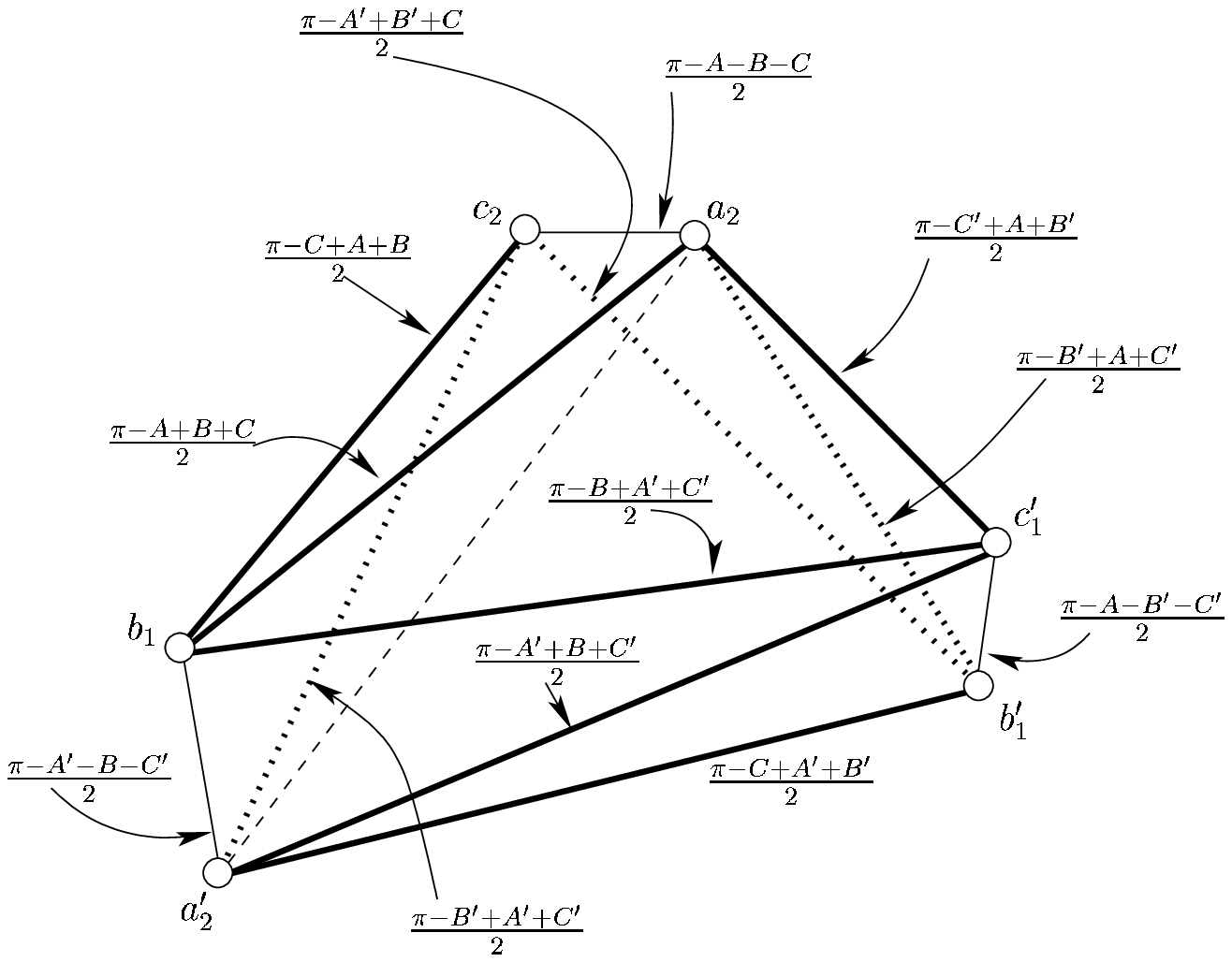}{f:octahedron}{Decomposition of
polyhedron $U$, Part II: the octahedron $O$}
\subsection{Triangulating the octahedron}\label{s:tri_oct}
The octahedron in Figure~\ref{f:octahedron} can be triangulated by drawing an edge
from one of its vertices to a non-adjacent vertex. There are three distinct ways
to do this. No matter which way is chosen, the dihedral angles of the resulting
four tetrahedra can be determined by a family of 7 linear equations and one
quadratic equation.

The following terminology will be used to describe the triangulation of an octahedron.
\begin{defn}\label{d:firepole}
The {\em firepole} of an octahedron is a segment joining one of its vertices
to a non-adjacent vertex.
\end{defn}
In Figure~\ref{f:octahedron} the firepole is the dashed line connecting vertices $a_2$
and $a_2^\prime$. At this point it is easier to view $O$  in the half-space
model, with the vertex $a_2$ at the point at infinity. This is depicted in
Figure~\ref{f:envelope}. The segments $\{a_2,c_1^\prime\}$, $\{a_2,b_1^\prime\}$,
$\{a_2,c_2\}$, $\{a_2,b_1\}$ become lines perpendicular
to the plane at infinity in the half-space model, and their projection is seen in
Figure~\ref{f:envelope} as vertices $v_a$, $v_b$, $v_c$, and $v_d$.
\bfig{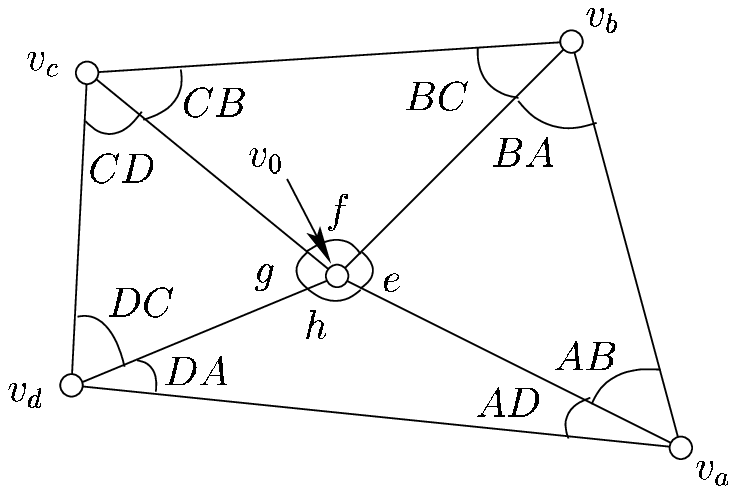}{f:envelope}{Octahedron $O$ in the half-space model}
The angles of the quadrilateral at each of these vertices are
$a$, $b$, $c$, and $d$, and they are equal to the
angles dihedral angles of $O$ at edges $\{v_0,v_a\}$, $\{v_0,v_b\}$, $\{v_0,v_c\}$, and
$\{v_0,v_d\}$, since the half-space model is conformal.
As seen in Figure~\ref{f:octahedron}, these angles are
\begin{equation}\label{e:abcd}
\begin{split}
a=\frac{\pi-C^\prime+A+B^\prime}{2};\;\;b=\frac{\pi-B^\prime+A+C^\prime}{2}\\
c=\frac{\pi-A-B-C}{2};\;\;d=\frac{\pi-A+B+C}{2}.
\end{split}
\end{equation}
The edges with dihedral angles $e$, $f$, $g$, $h$ are opposite to edges
$\{v_a,v_b\}$, $\{v_b,v_c\}$, $\{v_c,v_d\}$, $\{v_d,v_a\}$,
respectively, and the
angles at those edges have already been computed (see Figure~\ref{f:octahedron}). Using
the fact that opposite edges of an ideal tetrahedron have equal dihedral angles, we
have
\begin{equation}\label{e:efgh}
\begin{split}
e=\frac{\pi-A-B^\prime-C^\prime}{2};\;\;f=\frac{\pi-A^\prime+B^\prime+C}{2}\\
g=\frac{\pi-C+A+B}{2};\;\;h=\frac{\pi-B+A^\prime+C^\prime}{2}.
\end{split}
\end{equation}
The unknown angles have been denoted
as $AB$, $BA$, $BC$, $CB$, $CD$, $DC$, $DA$, and $AD$. These angles are subject to the
following linear constraints: \begin{equation}\label{e:linear_constraints} \begin{split}
AB+AD=a;\;\;&AB+BA+e=\pi\\ BC+BA=b;\;\;&BC+CB+f=\pi\\
CD+CB=c;\;\;&CD+DC+g=\pi\\
DA+DC=d;\;\;&DA+AD+h=\pi.
\end{split}
\end{equation}
The matrix expressing conditions \eqref{e:linear_constraints} has a
one-dimensional null space, so one more condition is needed to determine
the 8 unknown angles. Geometrically, this last condition insures that the four ideal
tetrahedra fit together. In other words, one can always fit together four
tetrahedra that satisfy \eqref{e:linear_constraints} as shown in
Figure~\ref{f:bad_envelope}, since the faces of all ideal tetrahedra are
ideal triangles, and, therefore, isometric to one another.
\bfig{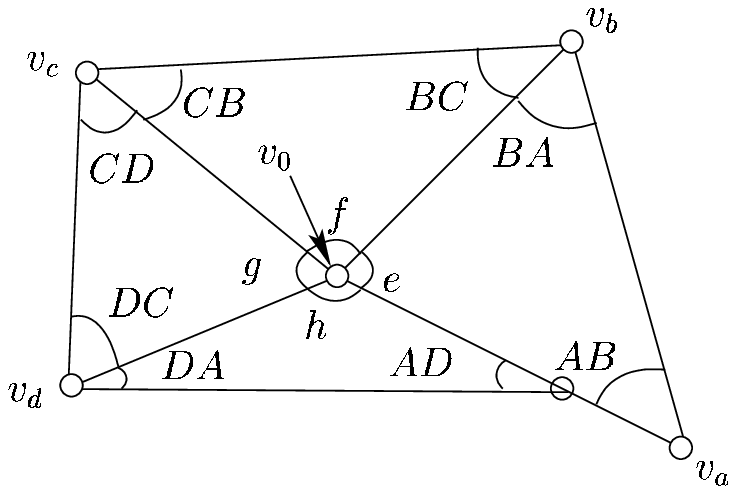}{f:bad_envelope}{Four tetrahedra that satisfy equations
\eqref{e:linear_constraints}}
In order avoid the situation in Figure~\ref{f:bad_envelope}, we need to
insure that once we fit the four tetrahedra together, the Euclidean length of
the segment $\{v_0,v_a\}$ does not change as we go around the quadrilateral
$\{v_a,v_b,v_c,v_d\}$. In particular, if we assume that the length of the
segment $\{v_0,v_a\}$ is equal to 1, then express it in terms of the other
lengths of the quadrilateral and equate the two quantities, we will get a
non-trivial condition that, together with equations
\eqref{e:linear_constraints}, will determine the unknown angles in
Figure~\ref{f:envelope}.

The equation we need is
\begin{equation}\label{e:holonomy}
\frac{\sin(AB)}{\sin(BA)}\frac{\sin(BC)}{\sin(CB)}
\frac{\sin(CD)}{\sin(DC)}\frac{\sin(DA)}{\sin(AD)}=1,
\end{equation}
and Figure~\ref{f:holonomy} suggests how it is obtained.
\bfig{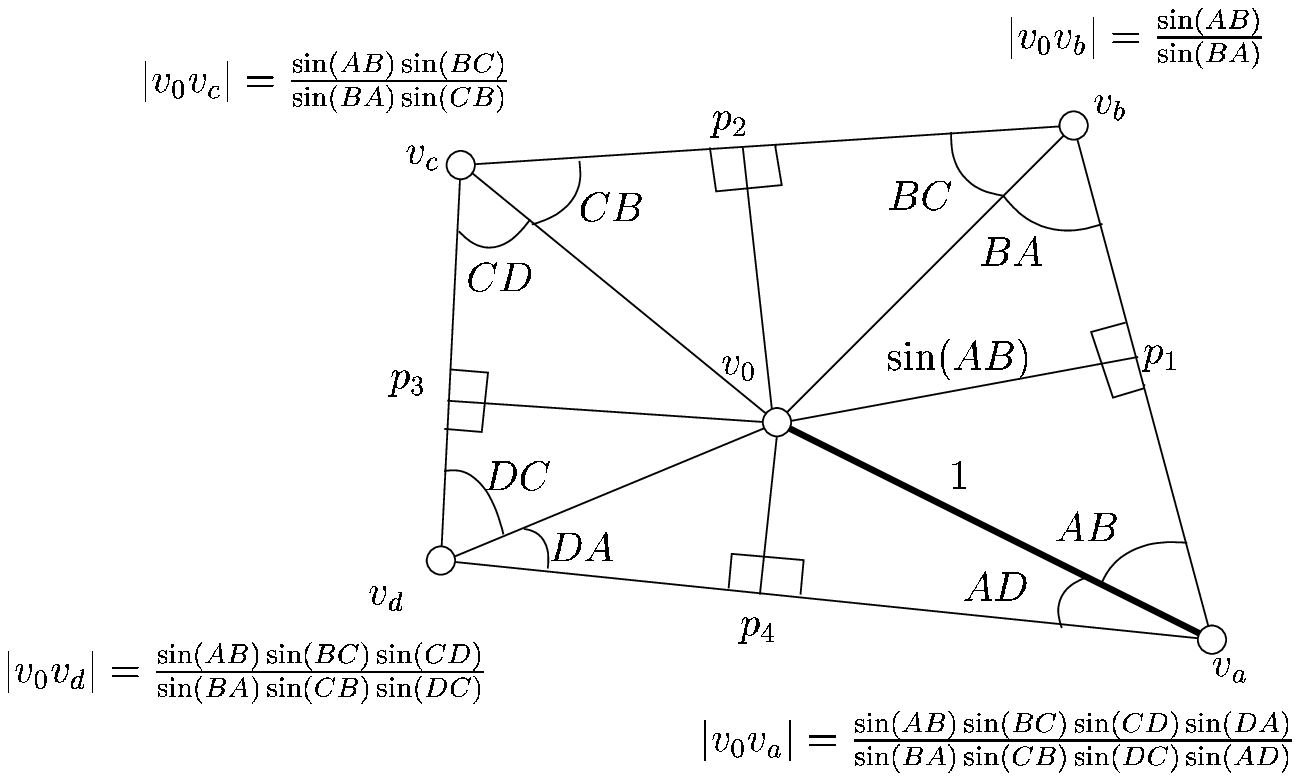}{f:holonomy}{Non-linear condition for insuring that four
tetrahedra fit together to make an octahedron}
Since scaling is an isometry in hyperbolic space, we may assume without loss of generality
that the Euclidean length of the segment $\{v_0,v_a\}$ is 1. Then by going
counter-clockwise around Figure~\ref{f:holonomy} and using basic trigonometry we
obtain
\begin{eqnarray*}
|v_0p_1|&=&\sin(AB)\\
|v_0v_b|&=&|v_0p_1|/\sin(BA)\\
        &=&\sin(AB)/\sin(BA),
\end{eqnarray*}
and so on, until
finally we arrive at the two equivalent expressions for $|v_0v_a|$ in
\eqref{e:holonomy}.

It is easier to solve the system of equations \eqref{e:linear_constraints} and
\eqref{e:holonomy} if we first come up with a solution in the one-dimensional space that
satisfies \eqref{e:linear_constraints} and then use \eqref{e:holonomy} to find the
remaining unknown. Let $(\overline{AB},\overline{BA},\overline{BC},\overline{CB},
\overline{CD},\overline{DC},\overline{DA},\overline{AD})$ be a solution to the system of
equations \eqref{e:linear_constraints}. Then there must be a $Z$ such
that
\begin{equation}\label{e:bars}
\begin{split}
\overline{AB}+Z=AB;\;\;&\overline{BA}-Z=BA\\
\overline{BC}+Z=BC;\;\;&\overline{CB}-Z=CB\\
\overline{CD}+Z=CD;\;\;&\overline{DC}-Z=DC\\
\overline{DA}+Z=DA;\;\;&\overline{AD}-Z=AD.
\end{split}
\end{equation}
No matter what the value of $Z$ is, the quantities on the right hand side of equations
\eqref{e:bars} still satisfy equations \eqref{e:linear_constraints}, since the sums in
those equations are unchanged when $Z$ is added to one summand and subtracted from the
other one.  After substituting \eqref{e:bars} into \eqref{e:holonomy},
letting $z=\exp(i Z)$ and
\begin{equation}\label{e:alphas_betas}
\begin{split}
\alpha_1=\exp(i \overline{AB});\;\;\beta_1=\exp(i \overline{BA})\\
\alpha_2=\exp(i \overline{BC});\;\;\beta_2=\exp(i \overline{CB})\\
\alpha_3=\exp(i \overline{CD});\;\;\beta_3=\exp(i \overline{DC})\\
\alpha_4=\exp(i \overline{DA});\;\;\beta_4=\exp(i \overline{AD}),
\end{split}
\end{equation}
we obtain
\begin{multline}\label{e:holonomy_exp}
0=\frac{1}{\alpha_1\alpha_2\alpha_3\alpha_4}-\beta_1\beta_2\beta_3\beta_4
+z^2(-\frac{\alpha_1}{\alpha_2\alpha_3\alpha_4}-\frac{\alpha_2}{\alpha_1\alpha_3\alpha_4}
-\frac{\alpha_3}{\alpha_1\alpha_2\alpha_4}-\frac{\alpha_4}{\alpha_1\alpha_2\alpha_3}\\
+\frac{\beta_1\beta_2\beta_3}{\beta_4}+\frac{\beta_1\beta_2\beta_4}{\beta_3}
+\frac{\beta_1\beta_3\beta_4}{\beta_2}+\frac{\beta_2\beta_3\beta_4}{\beta_1})\\
+z^4(\frac{\alpha_1\alpha_2}{\alpha_3\alpha_4}+\frac{\alpha_1\alpha_3}{\alpha_2\alpha_4}
+\frac{\alpha_2\alpha_3}{\alpha_1\alpha_4}+\frac{\alpha_1\alpha_4}{\alpha_2\alpha_3}
+\frac{\alpha_2\alpha_4}{\alpha_1\alpha_3}+\frac{\alpha_3\alpha_4}{\alpha_1\alpha_2}\\
-\frac{\beta_1\beta_2}{\beta_3\beta_4}-\frac{\beta_1\beta_3}{\beta_2\beta_4}
-\frac{\beta_2\beta_3}{\beta_1\beta_4}-\frac{\beta_1\beta_4}{\beta_2\beta_3}
-\frac{\beta_2\beta_4}{\beta_1\beta_3}-\frac{\beta_3\beta_4}{\beta_1\beta_2})\\
+z^6(\frac{\beta_1}{\beta_2\beta_3\beta_4}+\frac{\beta_2}{\beta_1\beta_3\beta_4}
+\frac{\beta_3}{\beta_1\beta_2\beta_4}+\frac{\beta_4}{\beta_1\beta_2\beta_3}\\
-\frac{\alpha_1\alpha_2\alpha_3}{\alpha_4}-\frac{\alpha_1\alpha_2\alpha_4}{\alpha_3}
-\frac{\alpha_1\alpha_3\alpha_4}{\alpha_2}-\frac{\alpha_2\alpha_3\alpha_4}{\alpha_1})
+z^8(\alpha_1\alpha_2\alpha_3\alpha_4-\frac{1}{\beta_1\beta_2\beta_3\beta_4}).
\end{multline}
By equations \eqref{e:linear_constraints} and \eqref{e:bars},
\begin{equation*}
(\overline{AB}+\overline{BC}+\overline{CD}+\overline{DA})+
(\overline{BA}+\overline{CB}+\overline{DC}+\overline{AD})=2\pi,
\end{equation*}
so that
\begin{equation*}
\alpha_1\alpha_2\alpha_3\alpha_4=\frac{1}{\beta_1\beta_2\beta_3\beta_4}.
\end{equation*}
This reduces \eqref{e:holonomy_exp} to a quadratic equaton in $z^2$.

Let $z_+$ ($z_-$) denote the solution of \eqref{e:holonomy_exp} corresponding to
adding (subtracting) the square root of the discriminant. Then $z_-$ gives the correct
values of the angles of the octahedron $O$ in Figure~\ref{f:octahedron}, while $z_+$ is of
great significance as well and will be discussed in \secref{ss:z_+}.

\subsection{Computation of the volume formula using the root $z_-$ of equation
\eqref{e:holonomy_exp}}
\label{ss:z_-}
At this point we can write down the volume of the octahedron $O$ (see
Figure~\ref{f:envelope}) using formula \eqref{e:allideal} for the
volume of an ideal hyperbolic tetrahedron:
\begin{eqnarray}\label{e:vol_octahedron_1}
V(O)&=&V(\{v_a,v_b,v_0,\infty\})+V(\{v_b,v_c,v_0,\infty\})+\nonumber\\
& &V(\{v_c,v_d,v_0,\infty\})+V(\{v_d,v_a,v_0,\infty\})\nonumber\\
&=&\Lob(AB)+\Lob(BA)+\Lob(e)+\Lob(BC)+\Lob(CB)+\Lob(f)\nonumber\\
& &\Lob(CD)+\Lob(DC)+\Lob(g)+\Lob(DA)+\Lob(AD)+\Lob(h).
\end{eqnarray}
Substituting \eqref{e:bars} and \eqref{e:efgh} into \eqref{e:vol_octahedron_1} and
setting $Z=\arg z_-$, we have
\begin{multline}\label{e:vol_octahedron_2}
V(O)=\Lob(\overline{AB}+\arg z_-)+\Lob(\overline{BA}-\arg z_-)+
\Lob(\overline{BC}+\arg z_-)\\+\Lob(\overline{CB}-\arg z_-)+
\Lob(\overline{CD}+\arg z_-)+\Lob(\overline{DC}-\arg z_-)\\+
\Lob(\overline{DA}+\arg z_-)+\Lob(\overline{AD}-\arg z_-)+
\Lob(\frac{\pi-A-B^\prime-C^\prime}{2})\\+\Lob(\frac{\pi-A^\prime+B^\prime+C}{2})+
\Lob(\frac{\pi-C+A+B}{2})+\Lob(\frac{\pi-B+A^\prime+C^\prime}{2}),
\end{multline}
where $\overline{AB}$, etc. are chosen as
\begin{equation}\label{e:ABs}
\begin{split}
&\overline{AB}=\frac{A+A^\prime+2B^\prime}{4};\;\;\overline{BA}=\frac{2\pi+A-A^\prime+2C^
\prime}{4}\\
&\overline{BC}=\frac{A+A^\prime-2B^\prime}{4};\;\;\overline{CB}=\frac{2\pi-A+A^\prime-2C}
{ 4}\\
&\overline{CD}=\frac{-A-A^\prime-2B}{4};\;\;\overline{DC}=\frac{2\pi-A+A^\prime+2C}{4}\\
&\overline{DA}=\frac{-A-A^\prime+2B}{4};\;\;\overline{AD}=\frac{2\pi+A-A^\prime-2C^\prime
} {4}. \end{split} \end{equation}
Using \eqref{e:vol_octahedron_2} along with the volumes of the three tetrahedra
and prism in Figure~\ref{f:U_decomp_tet_lab} and the fact that $\Lob$ is
odd and $\pi$-periodic gives \begin{eqnarray}\label{e:vol_U}
V(U)&=&\Lob(\overline{AB}+\arg z_-)+\Lob(\overline{BA}-\arg z_-)+\nonumber\\
& &\Lob(\overline{BC}+\arg z_-)+\Lob(\overline{CB}-\arg z_-)+\nonumber\\
& &\Lob(\overline{CD}+\arg z_-)+\Lob(\overline{DC}-\arg z_-)+\nonumber\\
& &\Lob(\overline{DA}+\arg z_-)+\Lob(\overline{AD}-\arg z_-)+\nonumber\\
& &+\Lob(A)+\Lob(A^\prime)+\Lob(B)+\Lob(B^\prime)+\Lob(C)+\Lob(C^\prime)\nonumber\\
& &+\Lob(\frac{\pi-A-B^\prime-C^\prime}{2})+\Lob(\frac{\pi+A^\prime-B-C^\prime}{2})
\nonumber\\
& &+\Lob(\frac{\pi+B^\prime-A-C^\prime}{2})+\Lob(\frac{\pi+C^\prime-A-B^\prime}{2})\nonumber\\
& &+\Lob(\frac{\pi+A-B-C}{2})+\Lob(\frac{\pi+C-A^\prime-B^\prime}{2})\nonumber\\
&
&+\Lob(\frac{\pi+B^\prime-A^\prime-C}{2})-\Lob(\frac{\pi+A^\prime+B^\prime+C}{2}
). \end{eqnarray} Finally, plugging \eqref{e:vol_U} into 
\eqref{e:unwrapped_candy} and using formula \eqref{e:3_ideal_c} for the volume
of an ideal prism gives \begin{multline}\label{e:greg_volume_222}
V(T)=\Lob(\overline{AB}+\arg z_-)+\Lob(\overline{BA}-\arg z_-)+
\Lob(\overline{BC}+\arg z_-)\\+\Lob(\overline{CB}-\arg z_-)+
\Lob(\overline{CD}+\arg z_-)+\Lob(\overline{DC}-\arg z_-)\\
+\Lob(\overline{DA}+\arg z_-)+\Lob(\overline{AD}-\arg z_-)+
\frac{1}{2}[\Lob(\frac{\pi+A-B-C}{2}\\-\Lob(\frac{\pi+B-A-C}{2})
-\Lob(\frac{\pi+C-A-B}{2})+\Lob(\frac{\pi+B^\prime-A^\prime-C}{2})\\
+\Lob(\frac{\pi+A+B+C}{2})+\Lob(\frac{\pi+C-A^\prime-B^\prime}{2})
+\Lob(\frac{\pi-A^\prime+B^\prime+C}{2})\\-\Lob(\frac{\pi+A^\prime+B^\prime+C}{2})
+\Lob(\frac{\pi+A^\prime-B-C^\prime}{2}-\Lob(\frac{\pi+A+B^\prime+C^\prime}{2})\\
-\Lob(\frac{\pi+A-B^\prime-C^\prime}{2})+\Lob(\frac{\pi+B^\prime-A-C^\prime}{2})
-\Lob(\frac{\pi-A^\prime-B+C^\prime}{2})\\+\Lob(\frac{\pi+A^\prime-B+C^\prime}{2
})
+\Lob(\frac{\pi+A^\prime+B+C^\prime}{2})+\Lob(\frac{\pi-A-B^\prime+C^\prime}{2})
], \end{multline} where the quantities with bars are given by \eqref{e:ABs}, and
$z_-$ is the solution of the quadratic equation \eqref{e:holonomy_exp} with the
negative square root. \subsection{Computation of the volume formula using the
root $z_+$ of equation \eqref{e:holonomy_exp}}
\label{ss:z_+}

When $z_-$ is replaced by $z_+$ in equation \eqref{e:greg_volume_222} one gets $-V(T)$
instead of $V(T)$. This surprising result has a very concrete geometrical explanation. The
main idea is that the octahedron $O$ has a dual octahedron $O^\prime$ associated with it,
and that $V(T)$ can be expressed in terms of either $O$ or $O^\prime$. It will be shown
below that the solution $z_+$ of  \eqref{e:holonomy_exp} solves the angles of a
triangulation of the octahedron $O^\prime$.

Recall formula \eqref{e:unwrapped_candy} in \secref{s:basic_setup}, which
indicates that we can construct a tetrahedron $T$ by subtracting 4 half-prisms
from the tetrahedron $U$ (see \figref{f:unwrapped_candy}). Let $H_P$ denote one
of the symmetric halves of a triangular prism $P$. So, for example, if $P$ is
the triangular prism $\{a,b,c,a^\prime,b^\prime,c^\prime\}$ in
\figref{f:prism_H}, $H_P=\{a,b,c,m_1,m_2,m_3\}$.
\afig{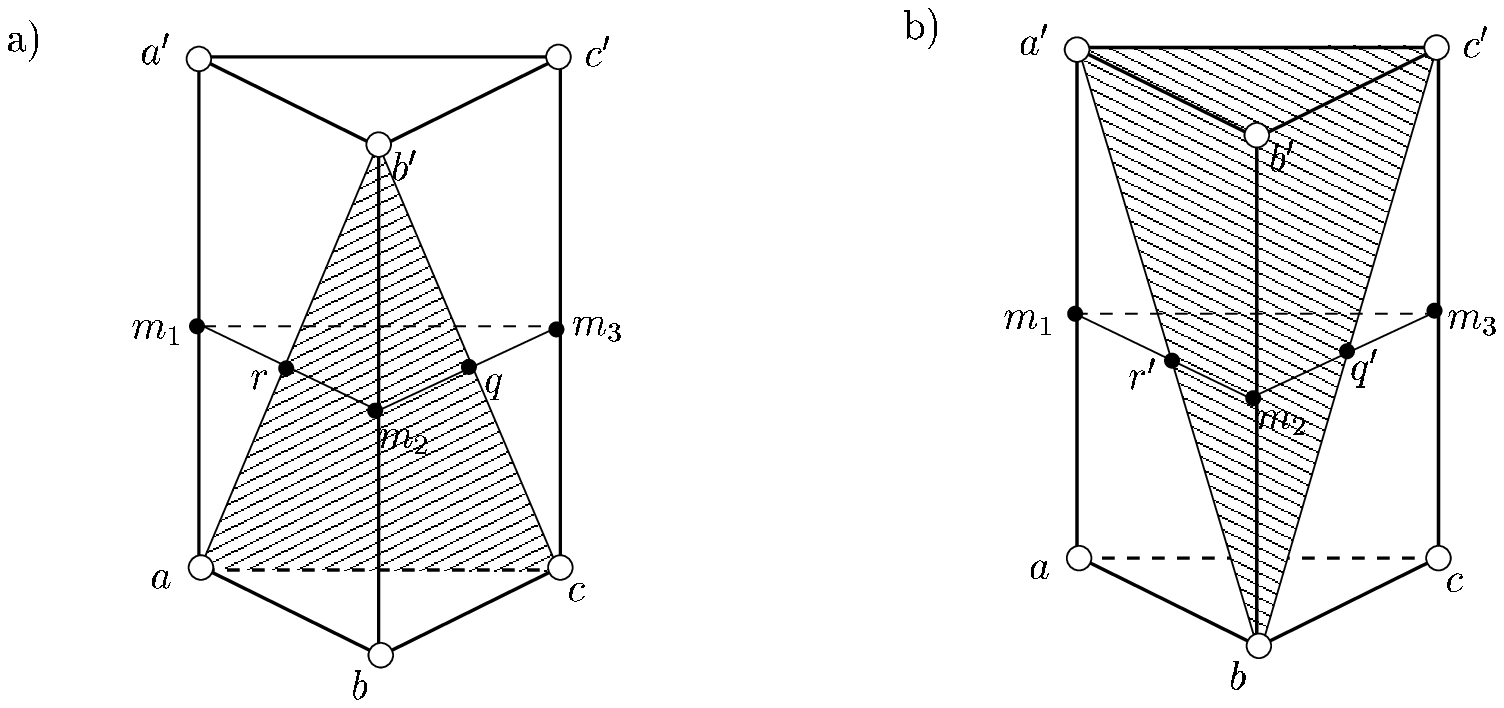}{f:prism_H}{Subdivisions of a triangular prism}
We will now show how to recast equation \eqref{e:unwrapped_candy} as a
decomposition. First, let us recall the construction of the octahedron $O$. As
described earlier, $O$ (see Figure~\ref{f:octahedron}) was constructed by
selecting one vertex from each edge of $U$ (see Figure~\ref{f:U_decomp_1}). It
follows that each of the four hexagonal faces of $U$ coincides with a face of
$O$. There is a one-to-one correspondence between the other four faces of $O$
and the four prisms \begin{eqnarray*}
P_1&=&\{c_1^\prime,a_1,b_1^\prime,c_2^\prime,a_2,b_2^\prime\}\\
P_2&=&\{a_1^\prime,b_2^\prime,c_1,a_2^\prime,b_1^\prime,c_2\}\\
P_3&=&\{a_2^\prime,b_1,c_2^\prime,a_1^\prime,b_2,c_1^\prime\}\\
P_4&=&\{a_2,b_2,c_2,a_1,b_1,c_1\}. \end{eqnarray*}
Formula \eqref{e:unwrapped_candy} can then be expressed as follows:
\begin{eqnarray}\label{e:vol_decomp_1}
T &=&U-H_{P_1}-H_{P_2}-H_{P_3}-H_{P_4}\nonumber\\
 &=&O+\nonumber\\
& &\{c_1^\prime,a_1,b_1^\prime,a_2\}+
\{a_1^\prime,b_2^\prime,c_1,a_2^\prime,b_1^\prime,c_2\}+
\{a_2^\prime,b_1,c_2^\prime,c_1^\prime\}+
\{a_2,b_2,c_2,b_1\}\nonumber\\
& &-H_{P_1}-H_{P_2}-H_{P_3}-H_{P_4}.
\end{eqnarray}
We now introduce $O^\prime$, the octahedron dual to $O$. With
Figure~\ref{f:U_decomp_1} in mind, one can visualize sliding the vertices
of $O$, $\{a_2,c_2,b_1,a_2^\prime,b_1^\prime,c_1^\prime\}$, along the respective edges
labeled as $A$, $B$, $C$, $A^\prime$, $B^\prime$, $C^\prime$ until they hit the
vertices at the end of these edges. We define the resulting octahedron as
$O^\prime$. In other words $O^\prime$ is the convex hull of the vertices
$\{a_1,c_1,b_2,a_1^\prime,b_2^\prime,c_2^\prime\}$. We can now express $U$ as
\begin{equation}\label{e:vol_U_halves}
U=O^\prime+\{c_1^\prime,b_1^\prime,b_2^\prime,c_2^\prime,a_1\}+
\{\emptyset\}+\{a_1^\prime, b_2,b_1, a_2^\prime,c_2^\prime\}+
\{a_1,a_2,c_2,c_1,b_2\},
\end{equation}
so that
\begin{multline}\label{e:vol_halves_2}
T=O^\prime+\{c_1^\prime,b_1^\prime,b_2^\prime,c_2^\prime,a_1\}+
\{\emptyset\}+\{a_1^\prime, b_2,b_1, a_2^\prime,c_2^\prime\}+
\{a_1,a_2,c_2,c_1,b_2\}\\
-H_{P_1}-H_{P_2}-H_{P_3}-H_{P_4}.
\end{multline}
We claim that
\begin{equation}\label{e:claim}
2T=O+O^\prime.
\end{equation}
To verify this claim, let us look back at
\figref{f:U_decomp_1} and how the last 8 terms of \eqref{e:vol_decomp_1}
are situated with respect to one another inside the polyhedron $U$. It is clear
that there are some partial cancellations between pairs of certain terms. For example, the polyhedra
$\{c_1^\prime,a_1,b_1^\prime,a_2\}$ and $H_{P_1}$ overlap. To see exactly what
happens, it is helpful to compare
$P_1=\{c_1^\prime,a_1,b_1^\prime,c_2^\prime,a_2,b_2^\prime\}$ to the prism
$\{a,b,c,a^\prime,b^\prime,c^\prime\}$ in \figref{f:prism_H}(a).
$\{c_1^\prime,a_1,b_1^\prime,a_2\}$ corresponds to $\{a,b,c,b^\prime\}$, and
$H_{P_1}$ corresponds to $\{a,b,c,m_1,m_2,m_3\}$. It follows immediately from
\figref{f:prism_H}(a) that \begin{equation}\label{e:div_tet}
\{a,b,c,b^\prime\}=\{a,b,c,r,m_2,q\}+\{r,m_2,q,b^\prime\}, \end{equation} and
\begin{equation}\label{e:div_pyr}
\{a,b,c,m_1,m_2,m_3\}=\{a,b,c,r,m_2,q\}+\{a,c,m_3,m_1,q,r\}. \end{equation}
Subtracting \eqref{e:div_pyr} from \eqref{e:div_tet} gives
\begin{equation}\label{e:div}
\{a,b,c,b^\prime\}-\{a,b,c,m_1,m_2,m_3\}=\{r,m_2,q,b^\prime\}-
\{a,c,m_3,m_1,q,r\}. \end{equation}
Now, the decomposition of $T$ in terms of $O^\prime$ in
\eqref{e:vol_halves_2} is obtained from the decomposition in
\eqref{e:vol_decomp_1} by sliding all the vertices of $O$ to the opposite ends
of the edges. So, for example, the tetrahedron
$\{c_1^\prime,a_1,b_1^\prime,a_2\}$ in \eqref{e:vol_decomp_1} gets replaced by
the pyramid $\{c_1^\prime,b_1^\prime,b_2^\prime,c_2^\prime,a_1\}$ in
\eqref{e:vol_halves_2}. One can see this in \figref{f:prism_H}, where the
tetrahedron $\{a,b,c,b^\prime\}$ corresponding to
$\{c_1^\prime,a_1,b_1^\prime,a_2\}$ is replaced by the pyramid
$\{a,c,c^\prime,a^\prime,b\}$ corresponding to
$\{c_1^\prime,b_1^\prime,b_2^\prime,c_2^\prime,a_1\}$ after the vertices $a$,
$c$, and $b^\prime$ are slid to $a^\prime$, $c^\prime$, and $b$.

As before, we wish to write the expression
$\{c_1^\prime,b_1^\prime,b_2^\prime,c_2^\prime,a_1\}-H_{P_1}$ in
\eqref{e:vol_halves_2} as a sum of disjoint simplices. \figref{f:prism_H}(b)
gives
\begin{equation*}
\{a,c,c^\prime,a^\prime,b\}=
\{a,b,c,m_1,r^\prime,m_2,q^\prime,m_3\}+\{a^\prime,c^\prime,m_1,m_3,
q^\prime,r^\prime\}
\end{equation*}
and
\begin{equation*}
\{a,b,c,m_1,m_2,m_3\}=
\{a,b,c,m_1,r^\prime,m_2,q^\prime,m_3\}+\{m_2,r^\prime,q^\prime,b\}.
\end{equation*}
Thus
\begin{equation}\label{e:evil_div}
\{a,c,c^\prime,a^\prime,b\}-\{a,b,c,m_1,m_2,m_3\}=
\{a^\prime,c^\prime,m_1,m_3,q^\prime,r^\prime\}-\{m_2,r^\prime,q^\prime,b\}.
\end{equation}
Comparing the right hand sides of \eqref{e:div} and \eqref{e:evil_div} we find
that they are negatives of one another. This is because $\{r,m_2,q,b^\prime\}$
and $\{a,c,m_3,m_1,q,r\}$ are isometric to $\{m_2,r^\prime,q^\prime,b\}$ and
$\{a^\prime,c^\prime,m_1,m_3,q^\prime,r^\prime\}$, respectively, by the symmetry
of the prism. It follows that the same relationship holds between
the corresponding terms in \eqref{e:vol_decomp_1} and \eqref{e:vol_halves_2}:
\begin{equation*}
\{c_1^\prime,a_1,b_1^\prime,a_2\}-H_{P_1}=-[\{c_1^\prime,b_1^\prime,b_2^\prime,
c_2^\prime,a_1\}-H_{P_1}].
\end{equation*}
It is easy to verify that the analagous relationship holds between the remaining
pairs of terms in \eqref{e:vol_decomp_1} and \eqref{e:vol_halves_2}. Therefore,
\begin{multline*}
-[\{c_1^\prime,a_1,b_1^\prime,a_2\}+
\{a_1^\prime,b_2^\prime,c_1,a_2^\prime,b_1^\prime,c_2\}+
\{a_2^\prime,b_1,c_2^\prime,c_1^\prime\}+
\{a_2,b_2,c_2,b_1\}\\-H_{P_1}-H_{P_2}-H_{P_3}-H_{P_4}]=\\
\{c_1^\prime,b_1^\prime,b_2^\prime,c_2^\prime,a_1\}+
\{a_1^\prime, b_2,b_1, a_2^\prime,c_2^\prime\}+\{a_1,a_2,c_2,c_1,b_2\}
-H_{P_1}-H_{P_2}-H_{P_3}-H_{P_4}.
\end{multline*}
Equation \eqref{e:claim} then
follows from adding \eqref{e:vol_decomp_1} and \eqref{e:vol_halves_2}.

Thus we have demonstrated that twice a hyperbolic tetrahedron is scissors congruent to two
octahedra, $O$ and $O^\prime$. This fact will be used extensively in the next
section to prove that the Regge symmetry is a scissors congruence.

All that remains to be shown is that by replacing $z_-$ with $z_+$ in
\eqref{e:vol_octahedron_2} one is swapping $V(O)$ for $-V(O^\prime)$. To see
this, let us examine the quantitative relationship between $O$ and $O^\prime$.
Recall Figure~\ref{f:prism_H}, where we saw that $\{a,b,c,b^\prime\}$ is isometric
to $\{c^\prime,b^\prime,a^\prime,b\}$. As stated earlier, we can visualize the
process of getting from Figure~\ref{f:prism_H}(a) to Figure~\ref{f:prism_H}(b)
as sliding the vertices $a$, $b^\prime$, $c$ along the edges $\{a,a^\prime\}$,
$\{b,b^\prime\}$, $\{c,c^\prime\}$, respectively. If we measure the change in
the angle $\theta$ between the oriented planes determined by $\{a,c,b^\prime\}$
and $\{a,b,b^\prime\}$ during this process, we find that $\theta$ changes to
$\pi-\theta$ in Figure~\ref{f:prism_H}(b). Extending this process to $O$ and
$O^\prime$, we see that dihedral angles of $O$ and $O^\prime$ are supplementary,
as shown in the Klein model drawings in Figure~\ref{f:good_evil}.
\afig{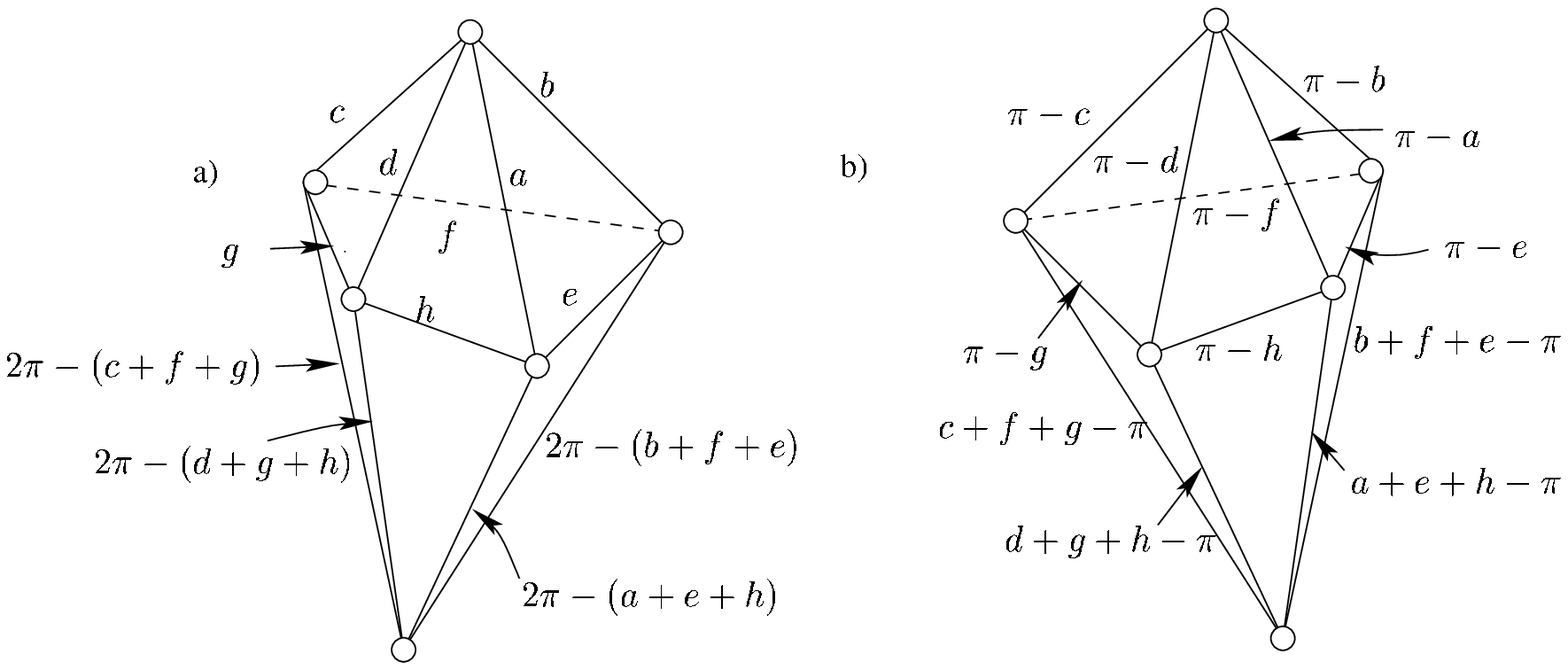}{f:good_evil}{a) The original octahedron $O$ with its
dihedral angles labeled\qua b) The dual octahedron $O^\prime$}
We can now apply the triangulation and computations described in
\secref{s:tri_oct} and \secref{ss:z_-} to $O^\prime$. This will yield a system
of linear equations similar to \eqref{e:bars} and a non-linear constraint like
\eqref{e:holonomy}. Choosing
\begin{equation*}
(-\overline{AB},-\overline{BC},-\overline{CD},-\overline{DA},
\pi-\overline{BA},\pi-\overline{CB},\pi-\overline{DC},\pi-\overline{AD})
\end{equation*}
as a solution to the system of linear equations results in a quadratic equation 
which turns out to be the same as \eqref{e:holonomy_exp}. It follows that
the second root of \eqref{e:holonomy_exp} solves for the unknown angles of the
triangulation of $O^\prime$. It is then easy to verify the following lemma.
\begin{lem}\label{l:zplusminus} Replacing $z_-$ with $z_+$ in
\eqref{e:vol_octahedron_2} results in the negative of the
volume of $O^\prime$.
\end{lem}
It follows that the root $z_+$ yields $-V(T)$ when substituted for $z_-$ in
\eqref{e:greg_volume_222}. We can use this fact to derive the following
expression for the volume of $T$:
\begin{multline}\label{e:clean_vol}
2V(T)=\Lob(\overline{AB}+\arg z_-)+\Lob(\overline{BA}-\arg z_-)+
\Lob(\overline{BC}+\arg z_-)\\+\Lob(\overline{CB}-\arg z_-)+
\Lob(\overline{CD}+\arg z_-)+\Lob(\overline{DC}-\arg z_-)\\+
\Lob(\overline{DA}+\arg z_-)+\Lob(\overline{AD}-\arg z_-)
+\Lob(-\overline{AB}+\arg z_+)\\-\Lob(\overline{BA}+\arg z_+)
+\Lob(-\overline{BC}+\arg z_+)-\Lob(\overline{CB}+\arg z_+)\\
+\Lob(-\overline{CD}+\arg z_+)-\Lob(\overline{DC}+\arg z_+)
+\Lob(-\overline{DA}+\arg z_+)\\-\Lob(\overline{AD}+\arg z_+).
\end{multline}

%
%
\section{Generating the Regge symmetries by scissors congruence}
\label{s:regge_proof}
Based on formulas \eqref{e:claim} and \eqref{e:clean_vol} we can
now construct a simple proof that $2T$ is scissors congruent to $2R(T)$, where
$R$ denotes any compositon of $R_a$, $R_b$, and $R_c$ as defined in
\eqref{e:R_a}, \eqref{e:R_b}, and \eqref{e:R_c}.

Recall Figure~\ref{f:envelope}, which shows a triangulation of $O$ in the half-space
model. In other words,
\begin{equation}\label{e:hull_octahedron}
O=\{v_a,v_b,v_0,\infty\}+\{v_b,v_c,v_0,\infty\}+\{v_c,v_d,v_0,\infty\}+\{v_d,v_a,v_0,
\infty\}.
\end{equation}
We will now subdivide each of the tetrahedra
$\{v_a,v_b,v_0,\infty\}$, $\{v_b,v_c,v_0,\infty\}$, \newline
$\{v_c,v_d,v_0,\infty\}$, and $\{v_d,v_a,v_0, \infty\}$ into three tetrahedra just
as Milnor did in his derivation of \eqref{e:allideal}. This is illustrated in Figure~\ref{f:double_tet}, where the dotted vertical line is a perpendicular
dropped from the vertex at the point of infinity, denoted by $\infty$, to the face opposite
to it, $\{a,b,c\}$.

\cfig{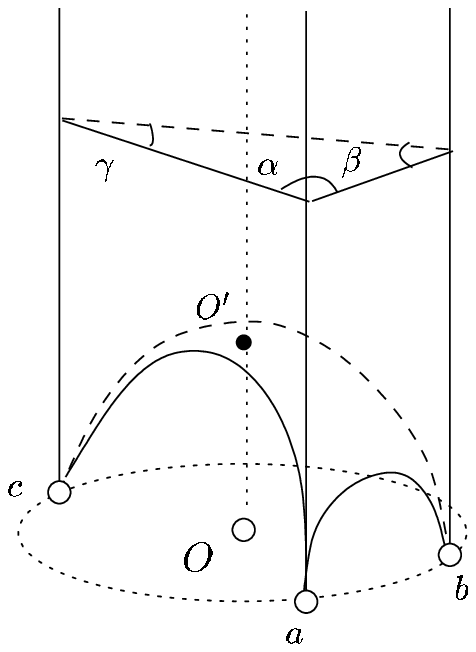}{f:double_tet}{An ideal hyperbolic tetrahedron in the half-space model}
This line must meet the
plane at infinity at the center of the hemisphere that determines $\{a,b,c\}$. In
other words, the projection of this configuration onto the plane at infinity looks
like Figure~\ref{f:milnor}, where the end of the perpendicular coincides with the
circumcenter of the triangle determined by $a$, $b$, and $c$.
\bfig{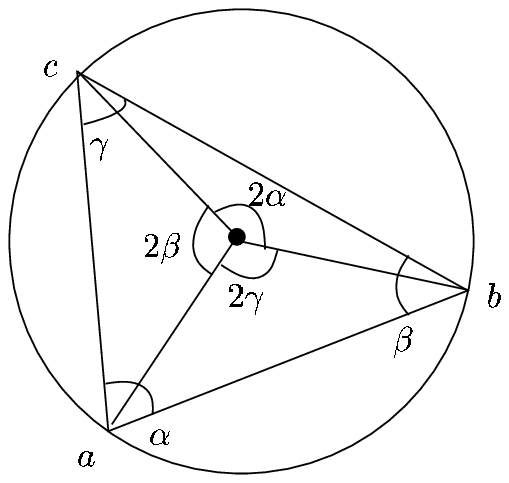}{f:milnor}{The projection of Figure~\ref{f:double_tet} onto the
plane at infinity}
Thus we see from \figref{f:double_tet} that $\{b,c,O^\prime,\infty\}$, $\{c,a,O^\prime,\infty\}$,
 and $\{a,b,O^\prime,\infty\}$ are $\L(\alpha)$, $\L(\beta)$, and $\L(\gamma)$, respectively.

We now
apply this construction to the 4 tetrahedra that triangulate the octahedron in
Figures~\ref{f:octahedron} and ~\ref{f:envelope}, ending up with 12 tetrahedra
each of which has 3 ideal vertices and 1 non-ideal vertex. The projection of
this construction onto the plane at infinity is shown in
Figure~\ref{f:milnored_envelope}.
 In that figure, the projections of the vertices $p_e$,
$p_f$, $p_g$, and $p_h$ are the respective circumcenters of the triangles
$v_av_bv_0$, $v_bv_cv_0$, $v_cv_dv_0$, and $v_dv_av_0 $. The actual positions of
the vertices $p_e$, $p_f$, $p_g$, and $p_h$ are on the planes determined by the
respective triangles, as are the dashed lines that represent the edges of the
new triangulation.
\afig{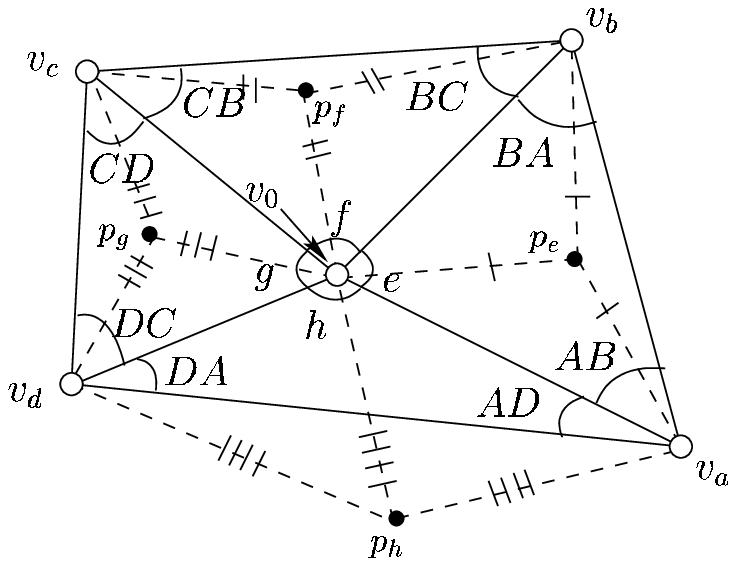}{f:milnored_envelope}{Octahedron $O$ in the
half-space model, triangulated according to the construction in
\figref{f:double_tet}}
\begin{rem} As Figure~\ref{f:milnored_envelope} indicates, the vertex
$p_h$ is outside of the triangle $v_dv_av_0$. As one might expect, formula
\eqref{e:allideal} still applies to the tetrahedron $\{v_d,v_a,v_0,\infty\}$.
Since the angle $h$ exceeds $\pi/2$, $\Lob(h)<0$. So in the formula
\begin{equation*} V(\{v_d,v_a,v_0,\infty\})=\Lob(AD)+\Lob(DA)+\Lob(h),
\end{equation*}
the first two terms correspond to the volumes of the tetrahedra $\{v_a,v_0,p_h,\infty\}$
and $\{v_0,v_d,p_h,\infty\}$, while the last term corresponds to the negative volume of
the tetrahedron $\{v_a,p_h,v_d,\infty\}$. Thus one can see how formula \eqref{e:allideal} still makes
geometric as well as analytic sense in the case of a tetrahedron such as $\{v_d,v_a,v_0,\infty\}$.
\end{rem}

Looking back at equations \eqref{e:vol_octahedron_1},
\eqref{e:claim}, and \figref{f:good_evil}, we see that the
terms $\Lob(e)$, $\Lob(f)$, $\Lob(g)$, $\Lob(h)$, which correspond to the
tetrahedra $\{v_a,v_b,p_e,\infty\}$, $\{v_b,v_c,p_f,\infty\}$,
$\{v_c,v_d,p_g,\infty\}$, $\{v_d,v_a,p_h, \infty\}$ in
Figure~\ref{f:milnored_envelope}, cancel in the formula for the volume of $T$.
Therefore, in our scissors congruence proof we will only be concerned with the
tetrahedra
\begin{equation*}
\begin{array}{cccc}
\{v_0,v_a,p_e,\infty\}&\{v_b,v_0,p_e,\infty\}&\{v_0,v_b,p_f,\infty\}&
\{v_c,v_0,p_f,\infty\}\\
\{v_0,v_c,p_g,\infty\}&\{v_d,v_0,p_g,\infty\}&\{v_0,v_d,p_h,\infty\}&
\{v_a,v_0,p _h,\infty\}
\end{array}
\end{equation*}
 and their
conterparts in $O^\prime$. These 8 tetrahedra are shown in
Figure~\ref{f:labeled_envelope}. They correspond to the first 4 terms in formula
\eqref{e:clean_vol}. \afig{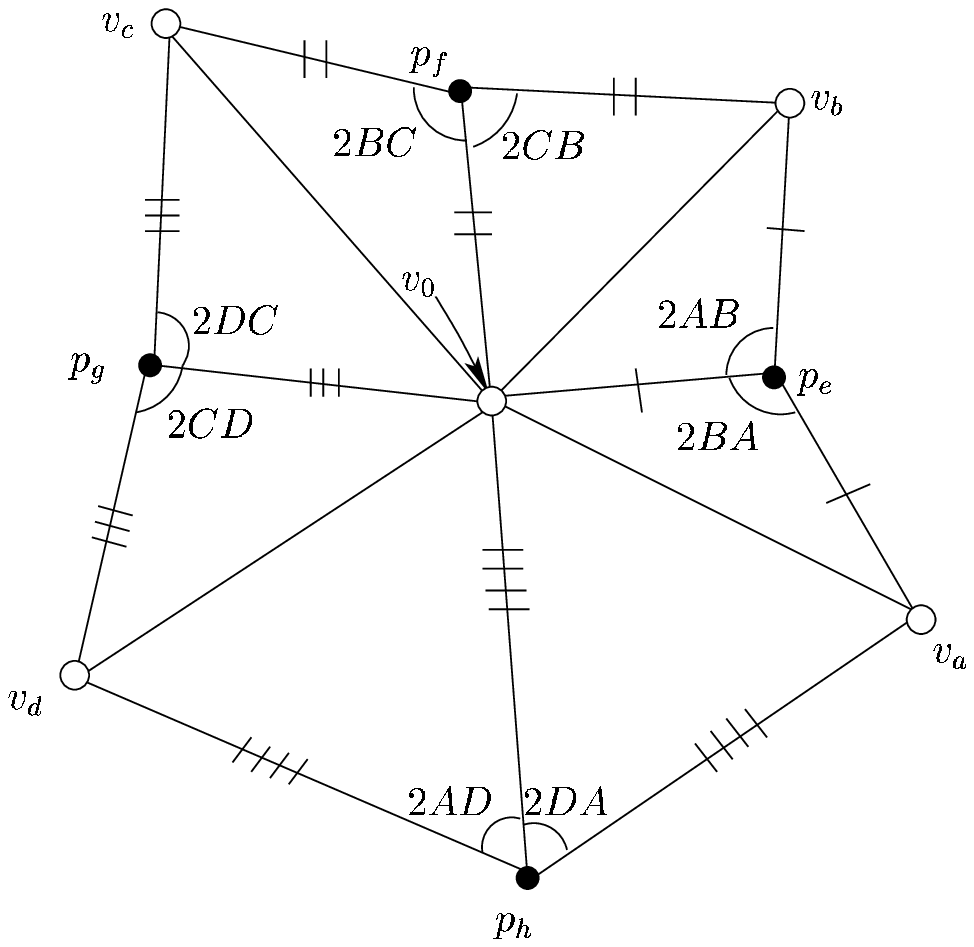}{f:labeled_envelope}{Tetrahedra
in the triangulation of the octahedron $O$ that correspond to the first 8 terms
in formula \eqref{e:clean_vol}} Triangulating $O^\prime$ in the same way as $O$
and applying the procedure described in \secref{s:tri_oct} and \secref{ss:z_-}
to find the unknown dihedral angles, we obtain the polyhedron shown in
Figure~\ref{f:labeled_evil_envelope}. The tetrahedra shown in that figure
correspond to the last 8 terms in formula \eqref{e:clean_vol}. Thus, up to a
multiple of $\pi$, $AB^\prime=-\overline{AB}+\arg z_+$,
$BA^\prime=-\overline{BA}-\arg z_+$, etc., where $\overline{AB}$,
$\overline{BA}$,... are given in \eqref{e:ABs} and $z_+$ is the positive square
root solution of \eqref{e:holonomy_exp}.
\afig{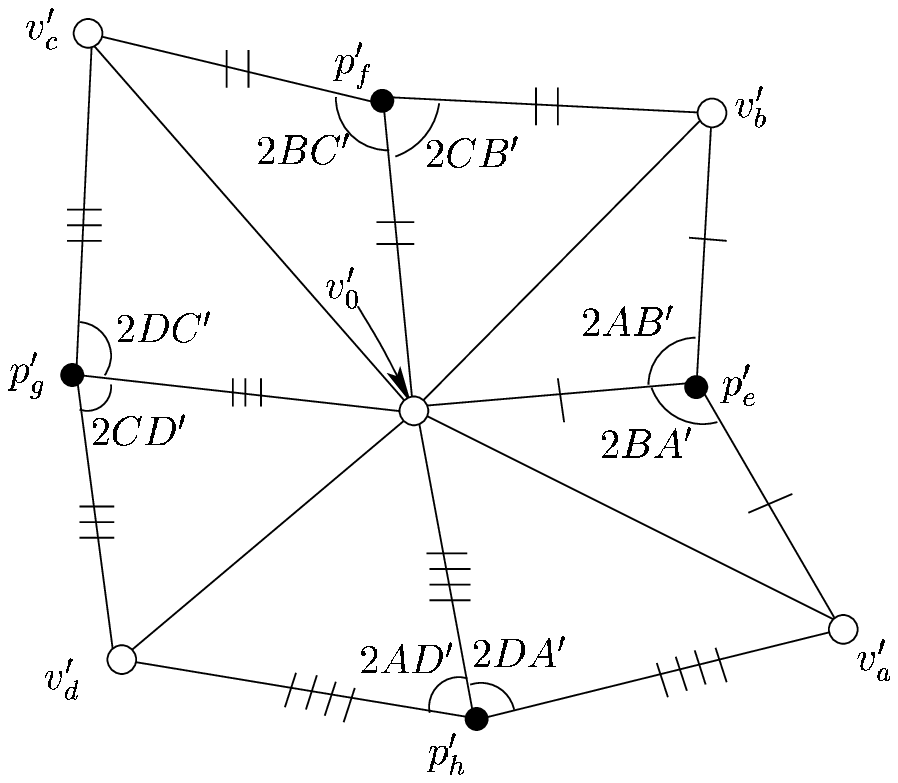}{f:labeled_evil_envelope}{Tetrahedra in the
triangulation of the octahedron $O'$ that correspond to the last 8 terms in
formula \eqref{e:clean_vol}} We have found that certain permutations of the
tetrahedra in Figure~\ref{f:labeled_envelope} (and corresponding permutations of
the tetrahedra in Figure~\ref{f:labeled_evil_envelope}) give us new polyhedra
that correspond to $R_b(T)$ and $R_c(T)$. In other words, given that $2T$ is
scissors congruent to $O+O^\prime$, we have found that $2R_b(T)$ is scissors
congruent to $P(O)+P(O^\prime)$, where $P(O)$ and $P(O^\prime)$ are the
octahedra obtained by permuting some of the tetrahedra in the triangulations of
$O$ and $O^\prime$.

Before stating this fact and its proof formally, we address the fact that there are no
permutations of the tetrahedra that give us $R_a(T)$. This can be traced back to the
choice of the firepole (see Definition~\ref{d:firepole}) in triangulating $O$.
The firepole, as shown in Figure~\ref{f:octahedron}, connects the vertices $a_2$ and
$a_2^\prime$. From the position of these vertices in Figure~\ref{f:unwrapped_candy},
we see that this choice of the firepole  ``favors'' the pair of opposite edges
labeled as $A$ and $A^\prime$. Similarly, the Regge symmetry $R_a$ is singles
out the edge with dihedral angle $A$ and its opposite, as seen in
Definition~\ref{d:regge}. The other two choices of the firepole for $O$ are the
segments $\{b_1,b_1^\prime\}$ and $\{c_1^\prime,c_2\}$, whose preferred
pairs of edges are labeled as $B$, $B^\prime$ and $C$, $C^\prime$. The first of
these choices yields a triangulation of $O$ that admits permutations that
correspond to $R_a(T)$ and $R_c(T)$, while the second allows permutations that
give $R_a(T)$ and $R_b(T)$. In fact, no matter which of the $2^6$ possible
decomposition one uses to cut down to octahedron, it is always true that for
every choice of a firepole that prefers a certain pair of opposite edges one
cannot obtain the octahedron corresponding to exactly one of the Regge
symmetries $R_a$, $R_b$ or $R_c$  by simply permuting tetrahedra.


The formal statement of the results obtained is as follows.
\begin{thm}\label{thm:regge_scissors}
Let $T=T(A,B,C,A^\prime,B^\prime,C^\prime)$ be a hyperbolic tetrahedron, and let
$R_b(T)=T(s_b-A,B,s_b-C,s_b-A^\prime,B^\prime,s_b-C^\prime)$ denote the action on $T$ of
one of the Regge symmetries, as defined in Definition~\ref{d:regge}. Then $2T$ is scissors
congruent to $2R_b(T)$.
\end{thm}
\begin{proof}
Let $O_T$ denote the octahedron obtained from $T$ by the construction in
\secref{s:gregs_formula}, and let $O_T^\prime$ denote the corresponding dual tetrahedron.
As shown in \secref{ss:z_+}, $2T$ is scissors congruent to $O_T+O_T^\prime$.
We apply the construction described in \secref{s:gregs_formula} to $R_b(T)$ by first
extending its edges to infinity and obtaining the octahedron $U_{R_b(T)}$ shown in
Figure~\ref{f:regge_b_unwrap}. Then we triangulate $U_{R_b(T)}$ so that $O_{R_b(T)}$ has
as its vertices the points $\{a_1,b_1^\prime,a_2^\prime,c_2^\prime,b_2,c_2\}$, while
$O_{R_b(T)}^\prime$ has vertices the
points $\{a_2,b_2^\prime,a_1^\prime,c_1^\prime,b_1,c_1\}$. $O_{R_b(T)}$ is depicted in
Figure~\ref{f:regge_b_oct}.
\bfig{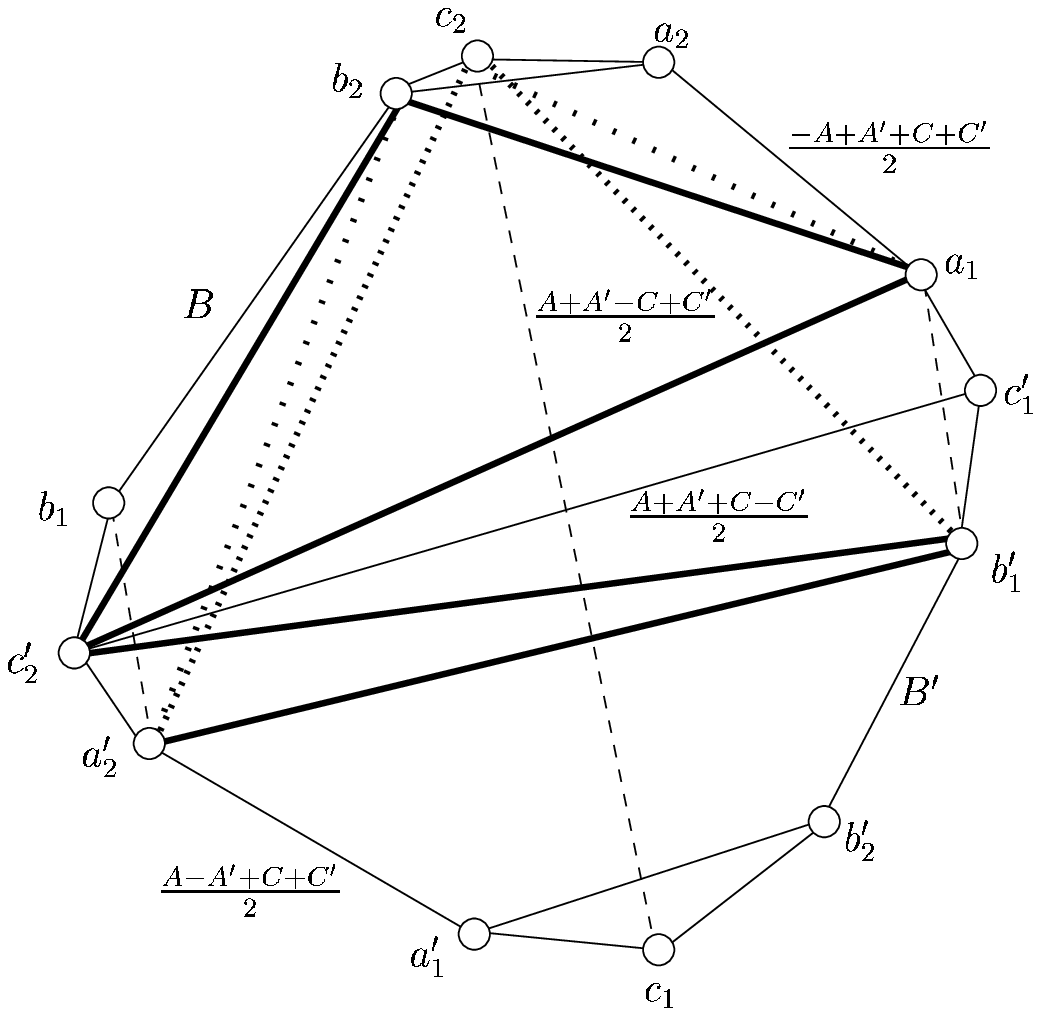}{f:regge_b_unwrap}{Polyhedron
$U_{R_b(T)}$}
\bfig{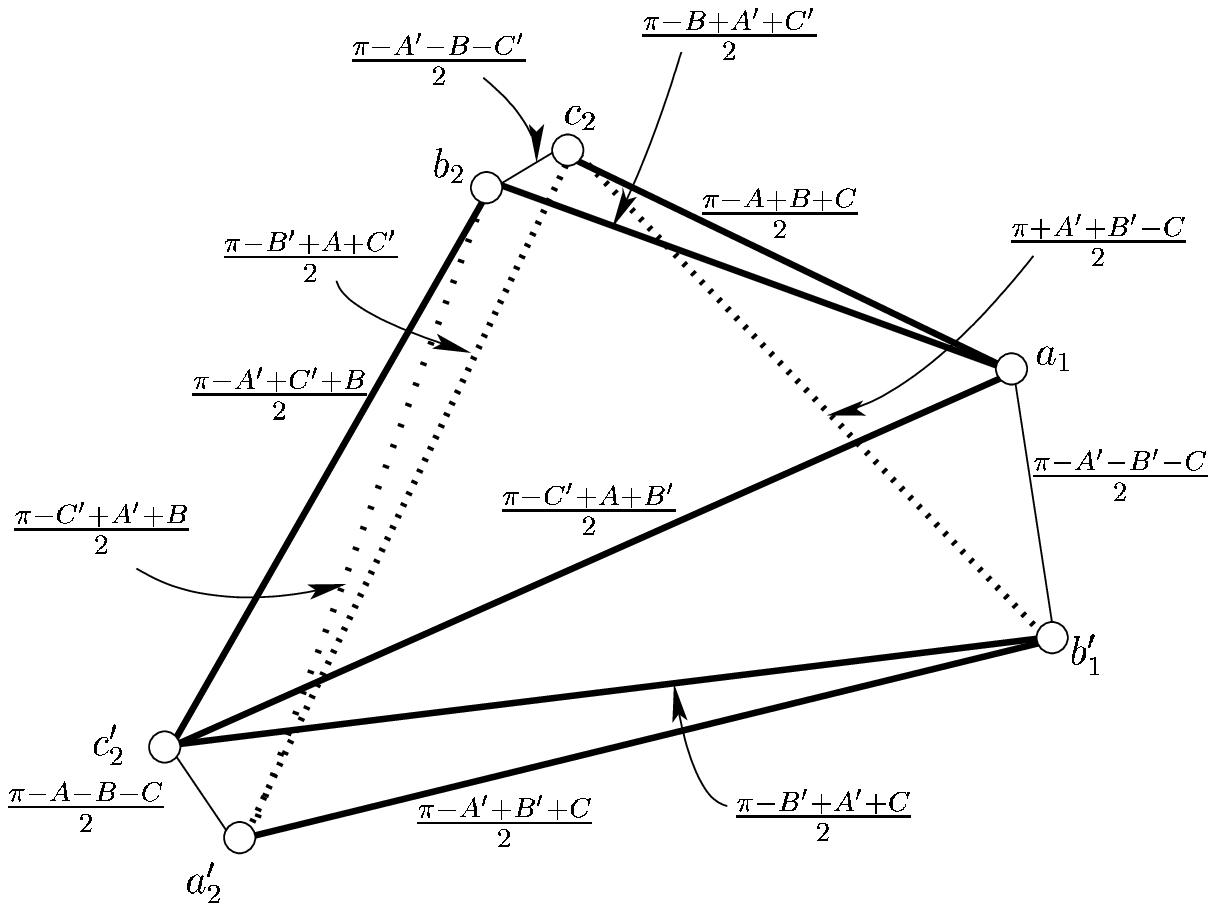}{f:regge_b_oct}{Octahedron $O_{R_b(T)}$}
Just as in the case of $T$, $2R_b(T)$ is scissors congruent to
$O_{R_b(T)}+O_{R_b(T)}^\prime$. What remains to be shown is that $O_T+O_T^\prime$ is
scissors congruent to $O_{R_b(T)}+O_{R_b(T)}^\prime$. This is done as follows.

Let $O$ be triangulated as shown in Figure~\ref{f:milnored_envelope}. We need not consider
the tetrahedra  $\{v_a,v_b,p_e,\infty\}$, $\{v_b,v_c,p_f,\infty\}$,
$\{v_c,v_d,p_g,\infty\}$, and $\{v_d,v_a,p_h,\infty\}$ as
they cancel when $O$ is added to $O^\prime$. Therefore, it
is sufficient to consider the polyhedron shown in Figure~\ref{f:labeled_envelope}. The
scissors congruence move consists of interchanging the tetrahedra $\{v_a,p_e,v_0,\infty\}$
and $\{v_0,v_c,p_g,\infty\}$, while leaving all the other tetrahedra in place. The
resulting polyhedron is shown in Figure~\ref{f:regge_b_envelope}, where the tetrahedra
that were moved are shaded.
\afig{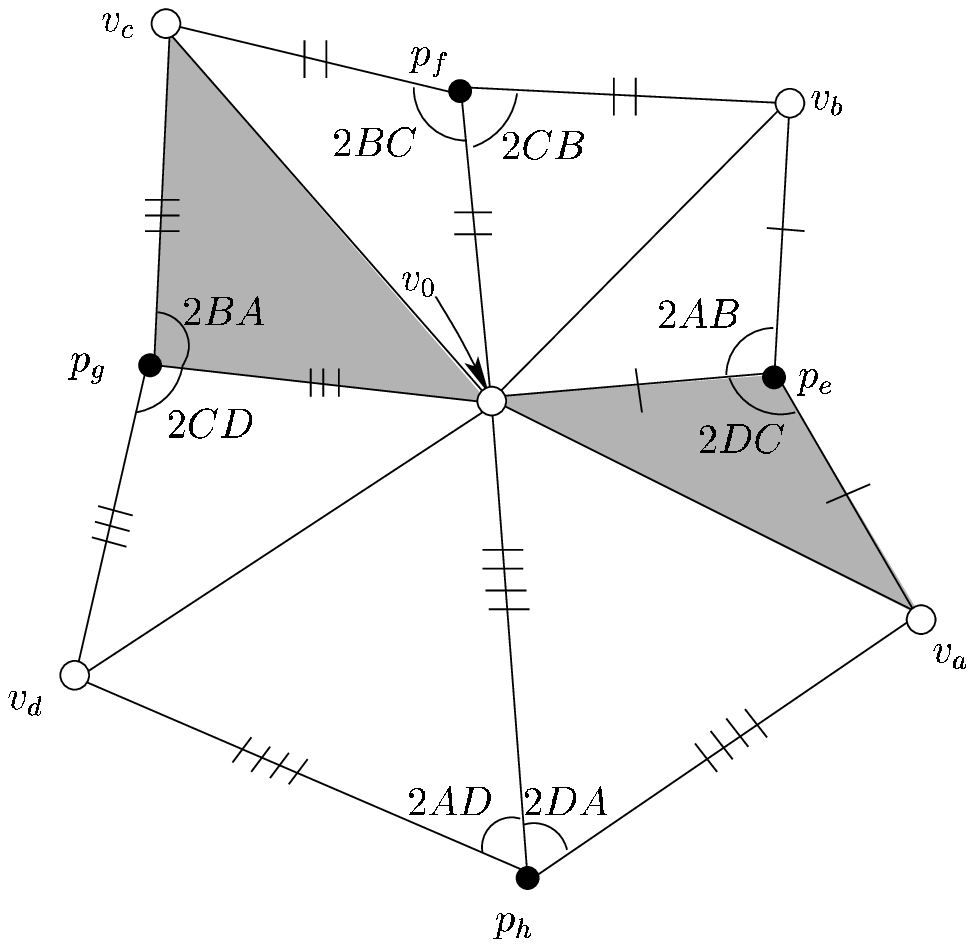}{f:regge_b_envelope}{Polyhedron of
Figure~\ref{f:labeled_envelope} with the interchanged tetrahedra shown shaded} 
The pieces in the resulting figure still fit together (i.e. the situation
depicted in Figure~\ref{f:bad_envelope} does not occur) because the new
polyhedron satisfies equation \eqref{e:holonomy}. In fact, the permutation
merely interchanges the terms $\sin(BA)$ and $\sin(DC)$ in \eqref{e:holonomy}.
As a final step, we translate the 8 tetrahedra making up the polyhedron of
Figure~\ref{f:regge_b_envelope} to obtain the mirror image of that polyhedron.
The result is shown in Figure~\ref{f:regge_b_envelope_mirror}.
\afig{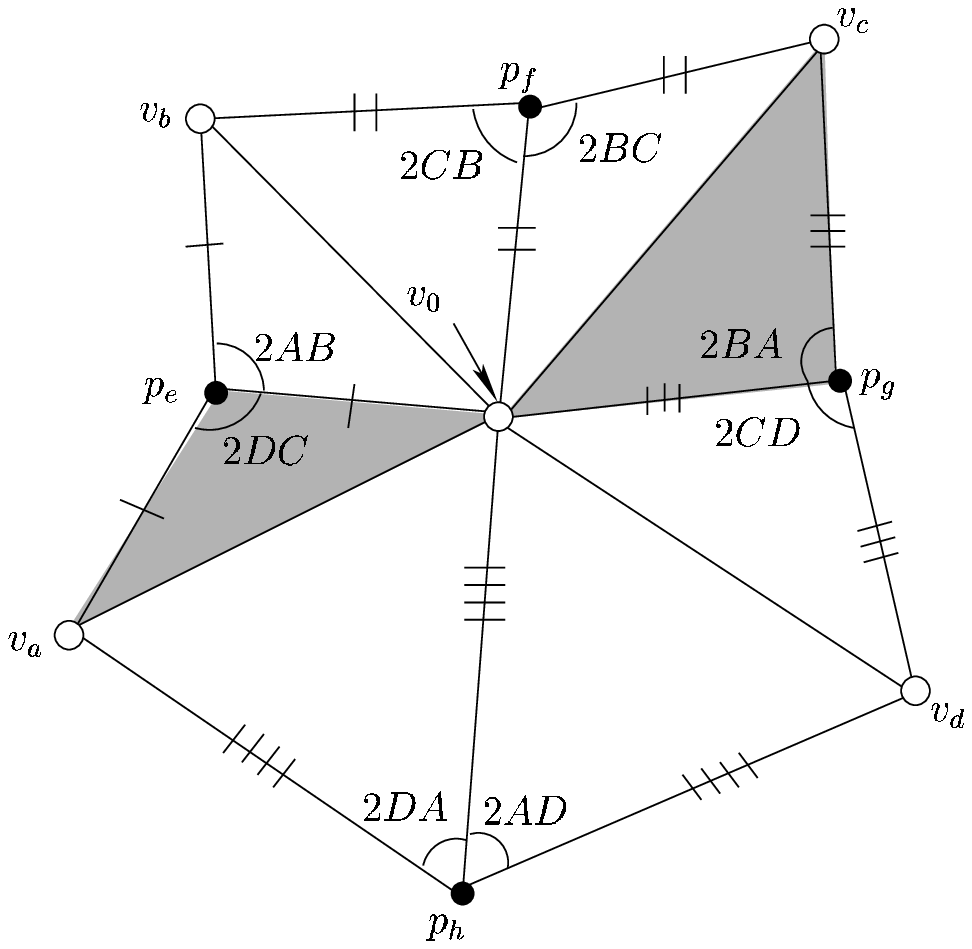}{f:regge_b_envelope_mirror}{Mirror image of
the polyhedron of Figure~\ref{f:regge_b_envelope}} After adding the tetrahedra
$\{v_b,v_a,p_e,\infty\}$, $\{v_c,v_b,p_f,\infty\}$, $\{v_d,v_c,p_g,\infty\}$,
and $\{v_a,v_d,p_h,\infty\}$ to the polyhedron in
Figure~\ref{f:regge_b_envelope_mirror} (these tetrahedra will subsequently
cancel with their counterparts in $O^\prime$) we obtain the octahedron $P(O)$
shown in \figref{f:regge_b_klein_mirror}.
The dihedral angles in that figure were obtained by adding up the dihedral
angles at the edges of the polyhedron in Figure~\ref{f:regge_b_envelope_mirror},
which were given in equations \eqref{e:bars}, \eqref{e:efgh} and
\eqref{e:ABs}. Clearly, $P(O)$ is isometric to the
octahedron $O_{R_b(T)}$ in Figure~\ref{f:regge_b_oct}.
\bfig{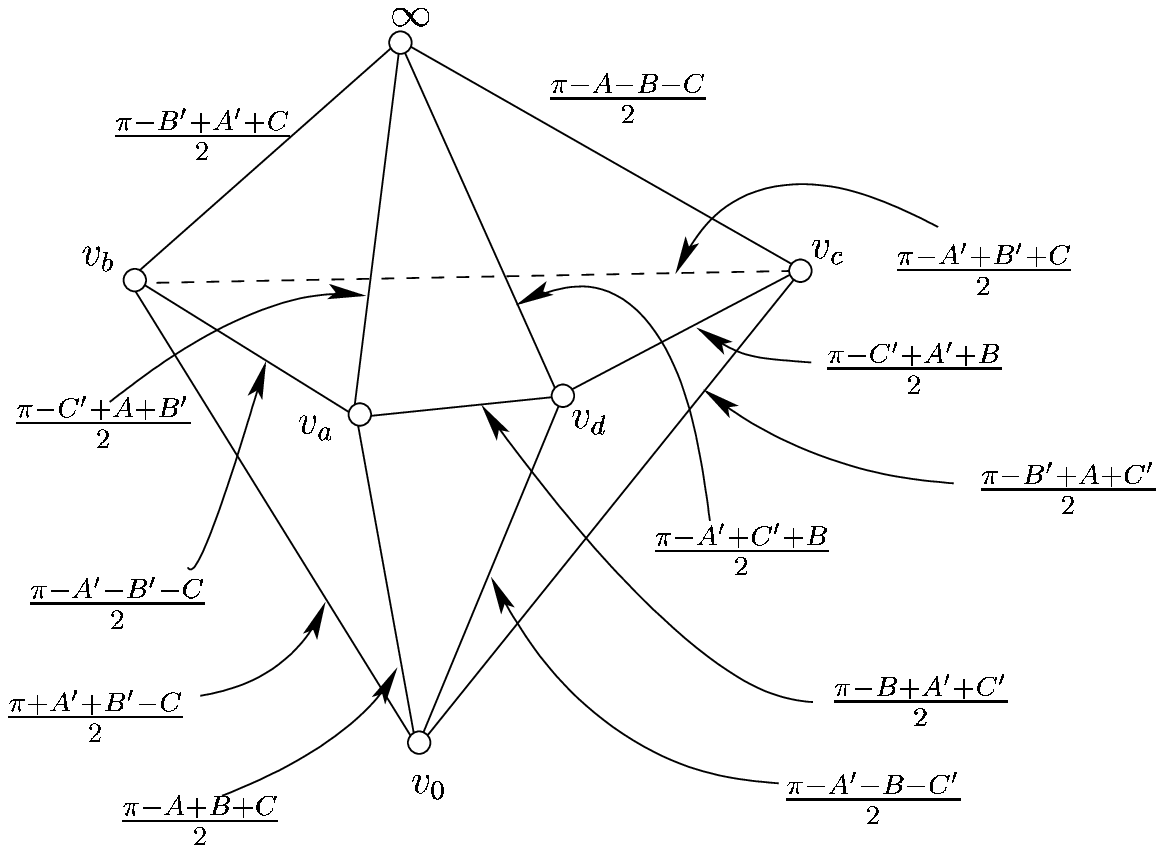}{f:regge_b_klein_mirror}{Octahedron $P(O)$ in the
Klein model} The discussion above applies verbatim to the dual octahedra. Since
their dihedral angles are dependent on the dihedral angles of the original
tetrahedra, we need to perform the same permutation on the tetrahedra making up
$O^\prime$ as we did on the tetrahedra making up $O$. In other words we
interchange the tetrahedra $\{v_a^\prime,p_e^\prime,v_0^\prime,\infty\}$ and
$\{v_0^\prime,v_c^\prime,p_g^\prime,\infty\}$ in
Figure~\ref{f:labeled_evil_envelope} and take the mirror image of the result. We
end up with the polyhedron shown in Figure~\ref{f:regge_b_evil_envelope_mirror}.
\afig{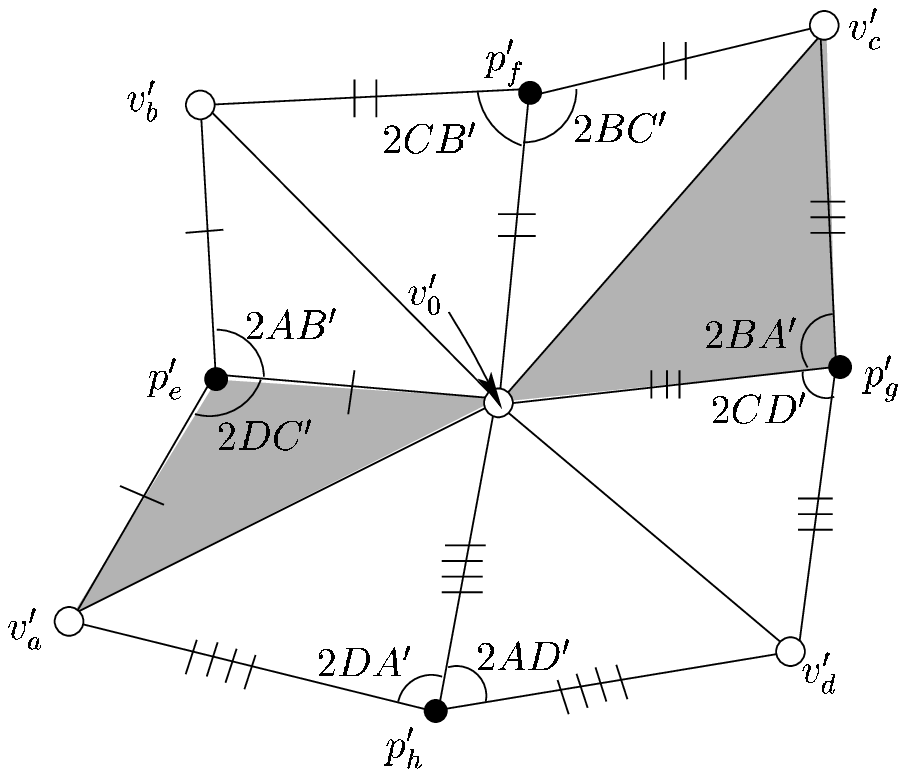}{f:regge_b_evil_envelope_mirror}{Mirror
image of the polyhedron of Figure~\ref{f:labeled_evil_envelope} taken after the
shaded tetrahedra have been interchanged}
Adding in the tetrahedra
$\{v_b^\prime,v_a^\prime,p_e^\prime,\infty\}$,
$\{v_c^\prime,v_b^\prime,p_f^\prime,\infty\}$,
$\{v_d^\prime,v_c^\prime,p_g^\prime,\infty\}$, and
$\{v_a^\prime,v_d^\prime,p_h^\prime,\infty\}$ gives the octahedron shown in
Figure~\ref{f:regge_b_evil_klein}.
\bfig{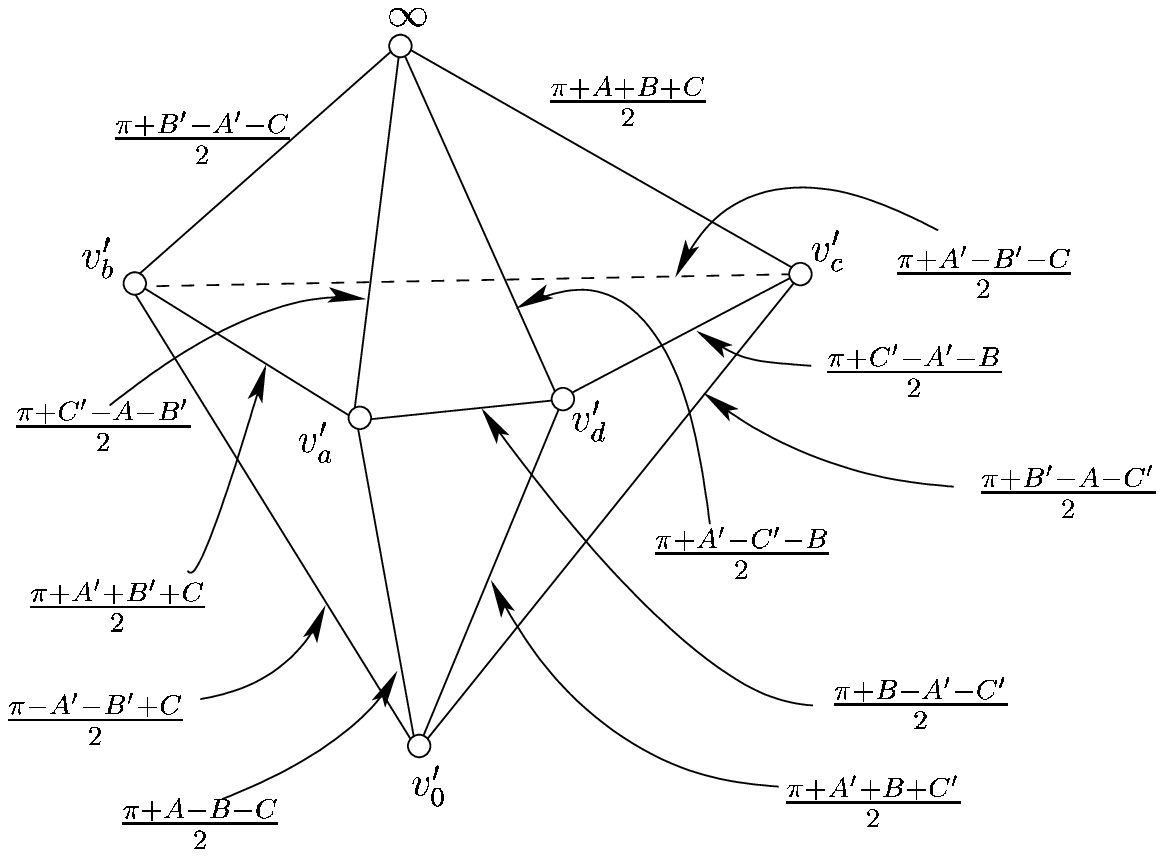}{f:regge_b_evil_klein}{Octahedron $P(O^\prime)$ in
the Klein model} Comparing the octahedra in Figure~\ref{f:regge_b_klein_mirror}
and Figure~\ref{f:regge_b_evil_klein}, we see that their dihedral angles add up
to $\pi$. Therefore, the octahedron in  Figure~\ref{f:regge_b_evil_klein} is
isometric to $O_{R_b(T)}^\prime$. Thus we have shown that $O+O^\prime$ is scissors
congruent to $O_{R_b(T)}+O_{R_b(T)}^\prime$ by a permutation. This completes the
proof. \end{proof} \begin{cor}
Let $T$ and $R_b(T)$ be hyperbolic tetrahedra as defined in \thmref{thm:regge_scissors}.
Then $T$ is scissors congruent to $R_b(T)$.
\end{cor}
\begin{proof}
The fact that we can ``divide by 2'' the construction that led to the proof of
\thmref{thm:regge_scissors} follows from Dupont's result of unique divisibility in
\cite{Dupont_book}. Dupont shows how an ideal hyperbolic tetrahedron can be divided into
two parts that are scissors congruent to one another. The fact that this 2-divisibility is unique means
that it is well-defined with respect to the different ways to divide a
tetrahedron into two scissors congruent sets of polyhedra. In other words, 
suppose we find that $T=A\coprod B$, where $A$ and $B$ are collections of
polyhedra such that $A$ is scissors congruent to $B$.
Now suppose that we can also express $T$ as
$T=A^\prime\coprod B^\prime$, where $A^\prime$ and $B^\prime$ are also scissors
congruent to one another. Then by unique 2-divisibility, $A$, $A^\prime$, $B$, and
$B^\prime$ are all scissors congruent to one another.

We saw in the proof of \thmref{thm:regge_scissors} that $2T$ is scissors congruent to
$\wt{O}+\wt{O}^\prime$, where $\wt{O}$ is the polyhedron in
Figure~\ref{f:labeled_envelope} and $\wt{O}^\prime$ is the polyhedron in
Figure~\ref{f:labeled_evil_envelope}. We now divide each of the 8 tetrahedra comprising
$\wt{O}$ into two scissors congruent (in fact, congruent) halves as shown in
Figure~\ref{f:div_by_2_envelope}.
\afig{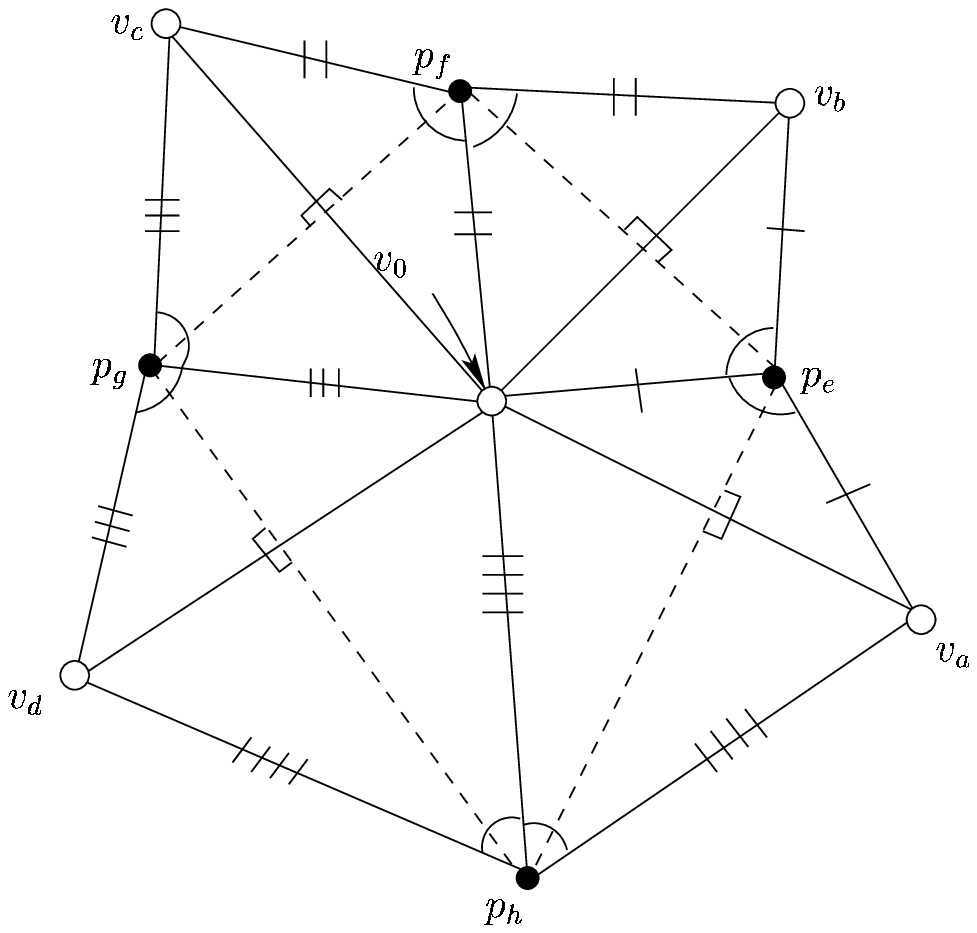}{f:div_by_2_envelope}{Dividing the polyhedron of
Figure~\ref{f:labeled_envelope} into two scissors congruent parts} The result 
is the polyhedron $\wt{O}/2$ shown in Figure~\ref{f:div_by_2_a}.
\bfig{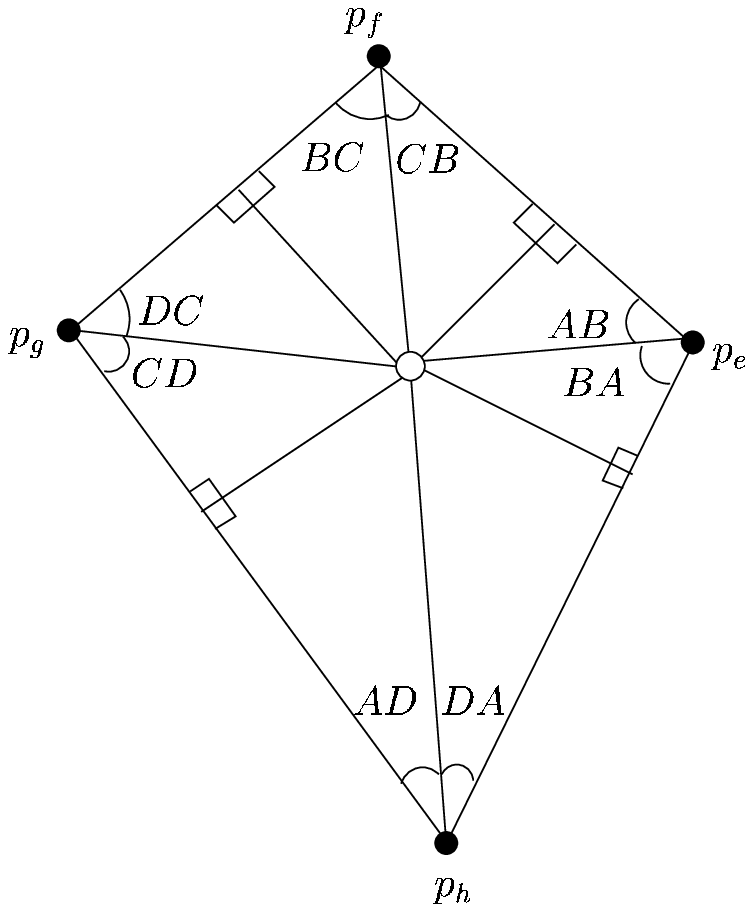}{f:div_by_2_a}{One of the two congruent halves comprising
the polyhedron of Figure~\ref{f:labeled_envelope}}
We perform the same division by 2 procedure on $\wt{O}^\prime$.

By Dupont's unique 2-divisibility result, $T$ is scissors congruent to
$\wt{O}/2+\wt{O}^\prime/2$, that is to the sum of the polyhedron in
Figure~\ref{f:div_by_2_a} and and its dual. Similarly, $R_b(T)$ is scissors congruent to
$\wt{O}_{R_b(T)}/2+\wt{O}_{R_b(T)}^\prime/2$. Permuting the indicated
tetrahedra in $\wt{O}/2$ results in $\wt{O}_{R_b(T)}/2$, shown in
Figure~\ref{f:div_by_2_regge_a}. Permuting the corresponding tetrahedra of
$\wt{O}^\prime/2$ results in $\wt{O}_{R_b(T)}^\prime/2$. Since
$\wt{O}/2+\wt{O}^\prime/2$ is scissors congruent to
$\wt{O}_{R_b(T)}/2+\wt{O}_{R_b(T)}^\prime/2$ by a permutation of tetrahedra, it
follows that $T$ is scissors congruent to $R_b(T)$.
\end{proof} \begin{cor} All the tetrahedra in the family generated by the Regge
symmetries are scissors congruent to one another.
\end{cor}
\begin{proof}
By \thmref{thm:regge_scissors}, the tetrahedron $T$ is scissors congruent to $R_b(T)$.
By applying the construction in the proof of that theorem to the rigid rotations of $T$ we
can prove that $T$ is scissors congruent to $R_a(T)$ and  $R_c(T)$. The corrollary then
follows from the fact that any two maps out of $R_a(T)$, $R_b(T)$ and  $R_c(T)$ generate
the group of Regge symmetries, so we can obtain any member of the family of Regge
symmetric tetrahedra by a sequence of procedures described in the proof of \thmref{thm:regge_scissors}.
\end{proof}
\bfig{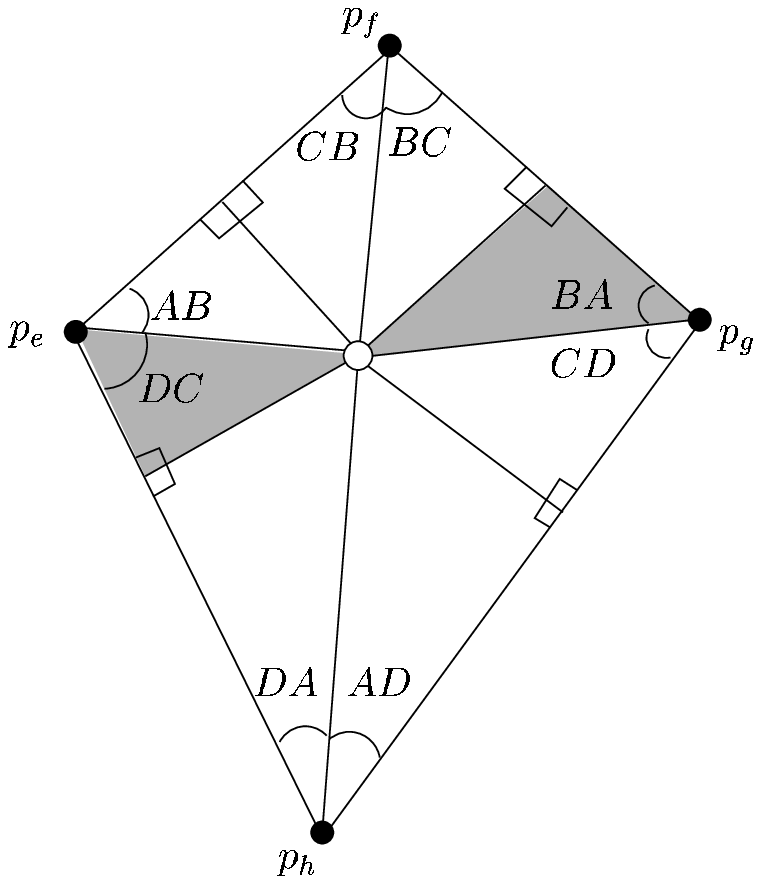}{f:div_by_2_regge_a}{One of the two congruent halves
comprising the polyhedron of Figure~\ref{f:regge_b_envelope_mirror}}
\section{Conclusion} \label{s:conclusion}
It would be natural at this point to try to extend the result of this paper to
spherical simplices. It is not even known with certainty
that the Regge symmetry is a scissors congruence in $\mathbb{S}^3$, since the
conjecture that the volume and Dehn invariant are sufficient invariants for
scissors congruence in spherical space is still open. Given the similarities
between spherical and hyperbolic spaces, it would certainly be reasonable to
hope that the Regge symmetries generate a family of scissors congruent
tetrahedra in $\mathbb{S}^3$.

Vinberg's result regarding the analytic continuation of the volume
function in \cite{Vinberg}, \S 3.1, tells us that Leibon's volume formulas
are valid for spherical tetrahedra up to multiplication by $\pm\sqrt{-1}$.
However, the difficulty in applying the construction described above to
spherical geometry lies in interpreting these formulas as geometric
decompositions. In particular, terms like $\Lob(\theta)$ can be viewed as
halves of the volume of bilaterally symmetric ideal tetrahedra in
hyperbolic space. But it is not clear what the analogue of an ideal
tetrahedron is in spherical geometry, and so the geometric significance of terms
like $\Lob(\theta)$ is difficult to see.

In Euclidean space, Roberts has already shown in \cite{Roberts} that the Regge
symmetry is a scissors congruence, but a constructive proof has not yet been
found. Adapting the above described construction to Euclidean space would
present a challenge that is in some sense more difficult than the spherical
case, since the volume formulas for Euclidean tetrahedra are so different from
those for hyperbolic tetrahedra. A successful adaptation of the construction to 
either spherical or Euclidean tetrahedra would certainly yield many insights
into simplices in both of these geometries.

\Addresses\recd
\end{document}